\documentclass[12pt]{article}

\usepackage{amsfonts}
\usepackage{amsmath}
\usepackage{amssymb}
\usepackage{graphics}
\usepackage{pictex}

\font\Bbb=msbm10
\def\X{\hbox{\Bbb X}}
\def\P{\hbox{\Bbb P}}
\def\C{\hbox{\Bbb C}}
\def\R{\hbox{\Bbb R}}
\def\Z{\hbox{\Bbb Z}}
\def\N{\hbox{\Bbb N}}

\newtheorem{theorem}{Theorem}
\newtheorem{proposition}[theorem]{Proposition}
\newtheorem{corollary}[theorem]{Corollary}
\newtheorem{lemma}[theorem]{Lemma}
\newtheorem{remark}{Remark}
\newtheorem{definition}{Definition}
\newtheorem{example}{Example}

\newenvironment{proof}{%
\smallbreak {\sc Proof.}}
{\smallbreak}

\begin{document}

\centerline{\bf\large SINGULARITIES OF NONCONFLUENT}

\vskip0.3cm

\centerline{\bf\large HYPERGEOMETRIC FUNCTIONS}

\vskip0.3cm

\centerline{\bf\large IN SEVERAL VARIABLES }

\vskip1cm

\centerline{Mikael Passare, Timur Sadykov and August Tsikh}

\vskip1cm

\noindent
{\small {\bf Abstract.}
The paper deals with singularities of nonconfluent hypergeometric
functions in several variables. Typically such a function is a
multi-valued analytic function with singularities along an algebraic
hypersurface. We describe such hypersurfaces in terms of
amoebas and the Newton polytopes of their defining polynomials.
In particular, we show that all $\mathcal{A}$-discriminantal
hypersurfaces (in the sense of Gelfand, Kapranov and Zelevinsky) have solid
amoebas, that is, amoebas with the minimal number of complement components.

\vskip0.3cm


\section{Introduction
\label{intr}}
There exist several approaches to the notion of hypergeometric series,
functions and systems of differential equations. In the present paper
we use the definition of these objects which was introduced by Horn
at the end of the 19th century~\cite{Horn}. His original definition
of a hypergeometric series is particularly attractive because of its
simplicity. A Laurent series in several variables is said to be
hypergeometric if the quotient of its two adjacent coefficients
depends rationally on the indices of summation.

In the present paper we study singularities of hypergeometric functions
which are defined by means of analytic continuation of hypergeometric
series. A hypergeometric series~$y(x)$ satisfies the so-called Horn
hypergeometric system
\begin{equation}
x_{i} P_i(\theta) y(x) = Q_i(\theta) y(x), \ i=1,\ldots ,n.
\label{horn}
\end{equation}
Here $P_{i}, Q_{i}$ are nonzero polynomials depending on the vector
differential operator $\theta=(\theta_{1},\ldots,\theta_{n}),$
$\theta_{i} = x_{i}\frac{\partial}{\partial x_{i}}.$
The nonconfluency of a hypergeometric series or the system~(\ref{horn})
means that the polynomials~$P_i$ and~$Q_i$ are of the same degree:
$$
{\rm deg }\, P_i = {\rm deg }\, Q_i, \,\,\, i=1,\ldots,n.
$$
These conditions can be expressed in terms of the Ore-Sato coefficient
of a hypergeometric series satisfying the system~(\ref{horn})
(see formulas~(\ref{hypcoeff}) and~(\ref{nonconfluency})).
Historically the Gauss hypergeometric differential equation was the
first one to be studied in detail due to the remarkable fact that
any linear homogeneous differential equation of order two with three
regular singularities can be reduced to it. The singularities of the
Gauss equation are $0,1,\infty.$ The generalized ordinary hypergeometric
differential equation which is a special case of the nonconfluent
system~(\ref{horn}) corresponding to $n=1$ also has three singular
points, namely $0,t,\infty,$ where~$t$ is the quotient of the
coefficients in the leading terms in the polynomials~$P_1$ and~$Q_1.$
Thus the singular set of an ordinary hypergeometric differential
equation is minimal in the following sense. There exist only two
circular domains, namely $\{0 < |x| < |t| \}$ and $\{ |t| < |x| < \infty \}$
in which any solution to the equation can be represented as a Laurent
series with the center at the origin (in the nonresonant case)
or as a linear combination of the products of Laurent series and
powers of~$\log x$ (in the resonant case).

It turns out that algebraic singularities of the system of partial
differential equations~(\ref{horn}) enjoy a multidimensional analogue
of this minimal property. It is convenient to formulate this property
in the language of amoebas which were introduced by Gelfand, Kapranov
and Zelevinsky in~\cite{GKZ}. The {\it amoeba} of an algebraic set
$\mathcal{R}=\{R(x)=0\}$ is defined to be its image under the mapping
${\rm Log } : (x_1,\ldots,x_n) \mapsto (\log|x_1|,\ldots,\log|x_n|).$
The complement of an amoeba consists of a finite number of convex
connected components which correspond to domains of convergence
of Laurent series expansions of single-valued functions with the
singularities on~$\mathcal{R}.$ The number of such components
cannot be smaller than the number of vertices of the Newton polytope
of the polynomial~$R(x).$ If these two numbers are equal then
we say that the amoeba is {\it solid.}
In Section~\ref{amoebasec} we prove the following theorem.

\vskip0.1cm
\noindent {\bf Theorem~\ref{rathypmin}}
{\it The singular hypersurface of any nonconfluent hypergeometric function
has a solid amoeba.}
\vskip0.1cm

A hypergeometric function satisfying the Gelfand-Kapranov-Zelevinsky
system of equations has singularities along the zero locus
of the corresponding $\mathcal{A}$-discriminant which is defined as
follows (see~\cite{Kapranov}).
Let~$\mathcal{A}$ be a finite subset of~$\Z^n$ and let~$f$
be a generic polynomial with the support~$\mathcal{A},$ i.e.,
$f=\sum\limits_{\alpha\in\mathcal{A}} c_{\alpha}x^{\alpha}.$
The corresponding $\mathcal{A}$-discriminant is defined to be
the polynomial in the coefficients~$c_{\alpha}$ which vanishes
whenever~$f$ together with all of its partial derivatives have a common
zero. Using Theorem~\ref{rathypmin} we arrive at the following corollary.

\vskip0.1cm
\noindent {\bf Corollary~\ref{adiscmin}}
{\it The zero set of any $\mathcal{A}$-discriminant has a solid amoeba. }
\vskip0.1cm

A geometric understanding of this latter result can be obtained from
the Horn-Kapranov uniformization theorem (see~\cite{Kapranov}), which states
that the logarithmic Gauss mapping on an $\mathcal{A}$-discriminantal
hypersurface is one-to-one. This implies that the normal directions of
the boundary of the corresponding amoeba are different at every
boundary point. In other words, two distinct tangent planes to the amoeba
boundary are never parallel. But if the amoeba complement were to
contain a bounded (convex) component there would have to be plenty of
distinct parallel tangent planes on the boundary.

Corollary~\ref{adiscmin} implies in particular that the amoeba
of the discriminant of a general algebraic equation is solid
(Corollary~\ref{discriminant}).

Let us also mention the following results in the paper.
Theorem~\ref{merorat} states that any meromorphic nonconfluent
hypergeometric function is rational.
In the last section we study the problem of describing the class of
rational hypergeometric functions. In the class of hypergeometric
functions satisfying the Gelfand-Kapranov-Zelevinsky system of equations
this problem was first considered in~\cite{CDD} and~\cite{CDS}.
Theorem~\ref{amoebathm}
gives a necessary condition for the Horn system to possess a rational
solution. The statement of Proposition~\ref{contigprop} emphasizes
the fact that only very few rational functions are hypergeometric.
The class of rational hypergeometric functions which is described
in this proposition consists of those which are contiguous to Bergman
kernels of complex ellipsoidal domains.

The proofs of the main results in the paper use the notions of the
support and the fan of a hypergeometric series, some facts from toric
geometry and the two-sided Abel lemma which is proved in
Section~\ref{meromorphicsec}.
Recall that the usual (one-sided) Abel lemma (see~\cite{GKZ0}
or~\cite{McDonald}) gives the following
relation between the domain of convergence
of a Puiseux series and its support (i.e., the set of summation).
\begin{lemma} {\rm (Abel's lemma for Puiseux series)}
Let $y(x)$ be a Puiseux series with a nonempty domain of convergence~$D.$
For any $x^{(0)}\in D$ and any cone~$C$ containing the convex hull of the
support of~$y(x)$ we have ${\rm Log }(x^{(0)}) - C^{\vee} \subset
{\rm Log }(D).$ Here~$C^{\vee}$ is the dual cone to~$C.$
\label{Abel}
\end{lemma}
The two-sided Abel lemma for hypergeometric Puiseux series states
that the domain~${\rm Log }\, (D)$ is itself contained in a suitable
translation of the cone~$-C^{\vee}.$

The authors are thankful to A.~Dickenstein for fruitful discussions.
Passare and Tsikh were supported by a grant from the Swedish Royal
Academy of Sciences. Sadykov and Tsikh were supported by the Russian
Foundation for Basic Research, grant 02-01-00167. Tsikh is grateful to
Max-Planck-Institut f\"ur Mathematik in Bonn for its hospitality and
financial support.


\section{Some basic notations and definitions
\label{resultantsec}}
To study the singularities of solutions to the Horn system~(\ref{horn})
we consider the characteristic variety of this system. Let~$\mathcal{D}$
denote the Weyl algebra of differential operators with polynomial
coefficients in~$n$ variables~\cite{Bjork1}. For any differential
operator $P \in \mathcal{D},$ $P=\sum_{|\alpha|\leq m} c_{\alpha}(x)
{\left( \frac{\partial}{\partial x}\right)}^{\alpha}$ its principal
symbol
$\sigma(P)(x,z)\in \C[x_{1},\ldots,x_{n},z_{1},\ldots,z_{n}]$
is defined by
$$
\sigma(P)(x,z)=\sum_{|\alpha|=m} c_{\alpha}(x)z^{\alpha}.
$$
We denote by~$G_{i}$ the differential operator
$x_{i} P_i(\theta) - Q_i(\theta)$ in the $i$th equation of the
Horn system~(\ref{horn}).
Let $\mathcal{M}= \mathcal{D}/ \sum_{i=1}^n \mathcal{D}G_{i}$ be the
left $\mathcal{D}$-module associated with the system~(\ref{horn}) and
let~$J\subset \mathcal{D}$ denote the left ideal generated by the
differential operators $G_{1},\ldots,G_{n}.$ By definition
(see~\cite{Bjork1}, Chapter~5, \S~2) the characteristic variety
${\rm char}(\mathcal{M})$ of the Horn system is given by
$$
{\rm char}(\mathcal{M})=\{ (x,z) \in \C^{2n} : \sigma(P)(x,z) =0,\,\,\,
{\rm for\ all\ } P\in J \}.
$$
We define the set~$U_{\mathcal{M}}\subset \C^n$ by
$$
U_{\mathcal{M}}=\{ x\in \C^n: \exists \, z\neq 0 {\rm\ such\ that\ }
(x,z)\in {\rm char}(\mathcal{M}) \}.
$$
It follows from Proposition~8.1.3 and Theorem~8.3.1 in~\cite{Hormander}
and Theorem~7.1 in Chapter~5 of~\cite{Bjork1} that a solution
to~(\ref{horn}) can only be singular on~$U_{\mathcal M}.$
Since any equation of the form $\sigma (P) (x,z) = 0$ is homogeneous
in~$z,$ it follows that~$U_{\mathcal{M}}$ is the image of
${\rm char}(\mathcal{M})$ under the projection of the direct product
$\C^n \times \P^{n-1} \rightarrow \C^n$ onto its first factor.
Using the main theorem of elimination theory (see \S~2C in~\cite{Mumford})
one can conclude that this image is an algebraic set, possibly the whole
of~$\C^n.$ In the latter case the singularities of a solution to the
Horn system are not necessarily algebraic. For instance, if every
differential operator~$G_i$ contains the factor $(\theta_1+\ldots+\theta_n)$
then any sufficiently smooth function depending on the quotients
$\frac{x_1}{x_n},\ldots,\frac{x_{n-1}}{x_n}$ is a solution to the
system~(\ref{horn}).

In the present paper we consider systems of the Horn type which
satisfy the condition $U_{\mathcal{M}}\neq\C^n.$ In this
case~$U_{\mathcal{M}}$ is a proper algebraic subset of~$\C^n.$
Its irreducible components of codimension greater than one are
removable as long as we are concerned with holomorphic solutions
to the Horn system. Thus the singular set of a solution to~(\ref{horn})
is algebraic and it is contained in the union of irreducible components
of codimension one. We denote this union by~$\mathcal{R}$ and call it
{\it the singular set of the Horn system.} Let~$R(x)$ be the defining
function of the set~$\mathcal{R},$ i.e.,
$$
\mathcal{R}= \{ R(x)=0 \}.
$$
The polynomial~$R(x)$ will be referred to as {\it the resultant of
the Horn system~(\ref{horn}).} To find a polynomial whose zero set
is~$\mathcal{R}$ is a difficult task which requires the full use of
elimination theory. There exists however a simple special case when
the set~$\mathcal{R}$ can be embedded into the zero set of some
polynomial which one can algorithmically compute.
Let $H_i(x,z)$ be the principal symbol of the differential operator~$G_i$
in the $i$th equation of the Horn system~(\ref{horn}). Since the
polynomials $H_{1},\ldots,H_{n}$ are homogeneous in $z_{1},\ldots,z_{n},$
they determine the classical resultant $R[H_1,\ldots,H_n]$
which is a polynomial in $x_1,\ldots,x_n$ (see~\cite{GKZ}, Chapter~13).
For the convenience of future reference we formulate the following
simple proposition.
\begin{proposition}
The singular set~$\mathcal{R}$ of the Horn system~(\ref{horn}) lies
in the zero set of the resultant $R[H_1,\ldots,H_n]$ of the principal
symbols of the operators in~(\ref{horn}).
\label{algebraicsing}
\end{proposition}
To prove this proposition it suffices to notice that for
$x^{(0)}\in U_{\mathcal M}$ the system of equations
$H_{1}(x^{(0)},z)=\ldots =H_{n}(x^{(0)},z)=0$ (considered as a system
of algebraic equations in $z_{1},\ldots,z_{n}$ whose coefficients
depend on~$x^{(0)}$) has a solution in~$\C^n \setminus \{0\}.$
This yields that the resultant of the homogeneous forms
$H_{1}(x,z),\ldots,H_{n}(x,z)$ with respect to the variables
$z_{1},\ldots,z_{n}$ vanishes at~$x^{(0)}$ (see~\cite{GKZ},
Chapter~13). Thus the singular locus of a solution to the Horn
system~(\ref{horn}) is contained in the zero set of the resultant
$R[H_1,\ldots,H_n].$ Notice that the vanishing of this resultant
at a point $x^{(0)}\in\C^n$ is equivalent to the condition that the
sequence of the principal symbols $\{ H_{i}(x^{(0)},z) \}_{i=1}^{n}$
is not regular in the polynomial ring $\C[z_{1},\ldots,z_{n}].$


\section{Puiseux series solutions to the Horn system and their supports
\label{supportsec}}

The Horn system~(\ref{horn}) as well as the Gelfand-Kapranov-Zelevinsky
system (see~\cite{GKZ0}) has the remarkable property that under some
natural assumptions there exists a basis in the space of its
holomorphic solutions consisting of (Puiseux) series with the center
at the origin (see~\cite{GKZ0} for the Gelfand-Kapranov-Zelevinsky
system and~\cite{Sadykov} for the Horn system). In this section we
introduce some terminology and present preliminary results which will
be used later for describing the singular set of the Horn system.

Suppose that a formal Puiseux series centered at the origin satisfies
the Horn system~(\ref{horn}). Such a series can be written as a linear
combination of formal shifted Laurent series, i.e., series of the form
\begin{equation}
y(x) =  x^{\gamma} \sum\limits_{s\in \Z^n} \varphi(s) x^s.
\label{series}
\end{equation}
Here $x^s=x_{1}^{s_{1}}\ldots x_{n}^{s_{n}}$, and the shift is 
determined by the
initial exponent
$\gamma =(\gamma_{1},\ldots,\gamma_{n}) \in \C^n, \,
{\rm Re}\,\gamma_{i} \in [0,1)$.
Suppose that the series~(\ref{series}) is a solution to~(\ref{horn}).
Computing the action of the operator $x_{i} P_i(\theta) - Q_i(\theta)$
on this series we arrive at the system of difference equations
\begin{equation}
\varphi(s+e_{i})Q_{i}(s+\gamma + e_{i}) =
\varphi(s)P_{i}(s+\gamma),
  \ i=1, \ldots ,n,
\label{difference}
\end{equation}
where~${\{e_{i}\}}_{i=1}^n$ is the standard basis of~$\Z^n.$
The system~(\ref{difference}) is equivalent to~(\ref{horn})
as long as we  are concerned with those solutions to the Horn system
which admit a series expansion of the form~(\ref{series}).

The system of difference equations~(\ref{difference}) is in general not
solvable without further restrictions on~$P_i,Q_i.$ Let~$R_i(s)$ denote
the rational function~$P_i(s)/Q_i(s+e_i),$ $i=1,\ldots,n.$ Increasing the
argument~$s$ in the~$i$th equation of~(\ref{difference}) by~$e_j$ and
multiplying the obtained equality by the~$j$th equation
of~(\ref{difference}), we arrive at the relation
$\varphi(s+e_i+e_j)/\varphi(s)=R_i(s+e_j)R_j(s).$ Similarly
$\varphi(s+e_i+e_j)/\varphi(s)=R_j(s+e_i)R_i(s).$
Thus the conditions
$R_{i}(s+e_{j})R_{j}(s) = R_{j}(s+e_{i})R_{i}(s),$ $i,j=1,\ldots,n$
are in general necessary for~(\ref{difference}) to be solvable.
Throughout this paper we assume that the polynomials~$P_i,Q_i$
defining the Horn system~(\ref{horn}) satisfy these relations and that
they are representable as products of linear factors.

The latter assumption together with the Ore-Sato theorem
(see~\cite{Sato} and~\cite{GGR}, \S~1.2) yields that the general solution
to the system of difference equations~(\ref{difference}) is of the
form
\begin{equation}
\varphi(s) = t_{1}^{s_{1}}\ldots t_{n}^{s_{n}} u(s)
\prod\limits_{i=1}^{p} \Gamma (\langle A_{i},s + \gamma \rangle - c_{i})
\phi(s).
\label{hypcoeff}
\end{equation}
Here $t_{i},c_{i}\in\C,$ $A_{i}=(A_{i1},\ldots,A_{in})\in\Z^n,$
$p\in\N_{0},$ $u(s)$ is a rational function whose numerator and
denominator are representable as products of linear factors and
$\phi(s)$ is an arbitrary periodic
function with the period~1 in each variable.
The fact that all the $\Gamma$-functions in~(\ref{hypcoeff}) are in
the numerator is unessential: using the identity
$\Gamma(z)\Gamma(1-z)=\pi / \sin \pi z$ and choosing the periodic
function~$\phi(s)$ in an appropriate way (see~\cite{Sadykov}), one can
move them into the denominator.
A formal series~(\ref{series}) with the coefficient~(\ref{hypcoeff})
is called a {\it formal solution} to the system~(\ref{horn}).
We will call any expression of the form~(\ref{hypcoeff}) the
{\it Ore-Sato coefficient} of a hypergeometric series (or of the
system~(\ref{horn})).
\begin{remark}
\rm
Conversely, the Ore-Sato coefficient~(\ref{hypcoeff}) defines
the system~(\ref{horn}) in the sense that for any $i=1,\ldots,n$
the quotient $\varphi(s+e_{i})/\varphi(s)$ equals
$P_{i}(s)/Q_{i}(s+e_{i}).$ For instance, the Ore-Sato
coefficient~(\ref{ex2soldiff}) in Example~\ref{example1} (see below)
defines the Horn system~(\ref{ex2horn}).
\label{definecoeff}
\end{remark}
The specific form of~(\ref{hypcoeff}) corresponds to our assumption
that the polynomials~$P_{i},Q_{i}$ can be represented as products
of linear factors. In general an Ore-Sato coefficient can include a
rational function which is not factorizable up to linear factors
(see~\cite{GGR}, \S~1.2). We may without loss of generality assume
that no linear factor in the rational function~$u(s)$ can be normalized
so that all of its coefficients become integers. Indeed, any linear
factor $a_{1}s_{1} + \ldots + a_{n}s_{n} + \lambda$ with $a_{i}\in\Z$
can be written in the form
$\Gamma(a_{1}s_{1} + \ldots + a_{n}s_{n} + \lambda + 1)/
\Gamma(a_{1}s_{1} + \ldots + a_{n}s_{n} + \lambda)$ and hence included
into the product of the $\Gamma$-functions in~(\ref{hypcoeff}).
Proposition~\ref{suppthm} (see below) yields that the other
linear factors of~$u(s)$ (such as $s_{1} + \pi s_{2}$)
are unessential as long as one is concerned with series solutions
to~(\ref{horn}). Throughout the paper we will assume that $u(s)\equiv 1.$

One can easily check that in terms of the parameters of the
Ore-Sato coefficient~$\varphi(s)$ the nonconfluency condition
${\rm deg}\, P_{i} = {\rm deg}\, Q_{i}$ can be written in the form
\begin{equation}
\sum_{i=1}^p A_{i} = 0.
\label{nonconfluency}
\end{equation}
Recall that in this paper we only deal with nonconfluent hypergeometric
series.

Any shifted Laurent series solution to~(\ref{horn}) (formal as well as
convergent) can be written in the form
\begin{equation}
y(x) =x^{\gamma}\sum\limits_{s\in S} \varphi(s) x^{s},
\label{solseries}
\end{equation}
where~$\varphi(s)$ is given by~(\ref{hypcoeff}) and~$S$ is a subset
of~$\Z^n$ on which $\varphi(s)\neq 0.$ The set~$S+\gamma$ will
be called the {\it support} of the series~(\ref{solseries}).
The support~$S+\gamma$ is
called {\it irreducible} if there exists no series solution
to~(\ref{horn}) supported in a proper nonempty subset of~$S+\gamma.$
A set~$S \subset \Z^n$ is said to be $\Z^n$-{\it connected} if any
two points of~$S$ can be connected by a polygonal line with unit
sides and vertices in~$S.$

Proposition~\ref{suppthm} (see below) describes all possible
supports of (formal) series solutions to~(\ref{horn}) and
Proposition~\ref{convergenceprop} allows one to find those of
them which have nonempty domains of convergence.
While looking for a solution to~(\ref{difference}) which is different
from zero on some subset~$S$
of~$\Z^n$  we will assume that the polynomials $P_{i}(s),Q_{i}(s),$
the set~$S$ and the vector~$\gamma$ satisfy the condition
\begin{equation}
|P_{i}(s+\gamma)|+|Q_{i}(s+\gamma +e_{i})| \neq 0,
\label{regular}
\end{equation}
for any~$s\in  S$ and for all $i=1, \ldots ,n.$ This assumption
eliminates the case when a solution to~(\ref{difference}) can
independently take arbitrary values at two adjacent points in the set~$S.$
The following statement (see~\cite{Sadykov}) gives necessary and
sufficient conditions for a solution to the system~(\ref{difference})
supported in some set~$S \subset \Z^n$ to exist.
\begin{proposition}
{\rm (Sadykov~\cite{Sadykov})}
For $S \subset \Z^n$ define
$$
S_{i}^{'} =\{s\in S: s+e_{i}\notin S \}, \,\,
S_{i}^{''} = \{s\notin S: s + e_{i} \in S\}, \,\,i=1, \ldots ,n.
$$
Suppose that the conditions~(\ref{regular}) are satisfied on~$S.$
Then there exists a solution to the system~(\ref{difference})  supported
in~$S$ if and only if the following conditions are fulfilled:
\begin{equation}
  P_{i}(s+\gamma) \vert_{S_{i}^{'}}    = 0, \,\,\,
  Q_{i}(s+\gamma +e_{i}) \vert_{S_{i}^{''}}    = 0, \,\, i=1, \ldots ,n,
\label{firstcond}
\end{equation}
\begin{equation}
  P_{i}(s+\gamma)
  \vert_{S\setminus   S_{i}^{'}}
\neq 0, \,\,\,
  Q_{i}(s+\gamma +e_{i}) \vert_{S}    \neq 0, \,\, i=1, \ldots ,n.
\label{secondcond}
\end{equation}
\label{suppthm}
\end{proposition}
By definition the union of the sets $S_{i}^{'}, S_{i}^{''},$
$i=1,\ldots,n$ is a discrete analogue of the boundary of the set~$S.$
Since the polynomials $P_i,Q_i$ are assumed to be representable as
products of linear factors, it follows from~(\ref{firstcond})
that~$S_{i}^{'}$ and~$S_{i}^{''}$ lie on hyperplanes. The
conditions~(\ref{secondcond}) yield that these hyperplanes bound
the set~$S.$ Thus we can formulate the following result.
\begin{proposition}
The convex hull of the support of a series solution to the Horn
system is a polyhedral set.
\label{polyhedral}
\end{proposition}
\begin{example}
\rm
Let us consider the following system of partial differential
equations of the Horn type
\begin{equation}
\left\{
\begin{array}{clcr}
x_{1} (\theta_{1} + \theta_{2})(\theta_{1} - 2) y(x) &=
(\theta_{1} - 1)(\theta_{1} - 4) y(x),\\
x_{2} (\theta_{1} + \theta_{2})(\theta_{2} - 3) y(x) &=
(\theta_{2} - 1)(\theta_{2} - 5) y(x).
\end{array}
\right.
\label{ex2horn}
\end{equation}
Assuming that~$y(x)$ admits a Laurent series expansion~(\ref{series})
with $\gamma=0,$ we arrive at the system of difference equations
\begin{equation}
\left\{
\begin{array}{clcr}
\varphi(s+e_{1}) s_{1}(s_1-3) =&
\varphi(s) (s_{1}+s_{2})(s_1-2),  \\
\varphi(s+e_{2}) s_{2}(s_2-4) =&
\varphi(s) (s_{1}+s_{2})(s_2-3).
\end{array}
\right.
\label{ex2diff}
\end{equation}
In accordance with the Ore-Sato theorem (see~\cite{Sato}
and~\cite{GGR}, \S~1.2) the general solution to the
system~(\ref{ex2diff}) is given by the function
\begin{equation}
\varphi(s) = (s_{1}-3)(s_{2}-4)
\frac
{\Gamma(s_{1}+s_{2})}
{\Gamma(s_{1})\Gamma(s_{2})}
\phi(s),
\label{ex2soldiff}
\end{equation}
where $\phi(s)$ is an arbitrary periodic function with the period~$1$
in~$s_{1}$ and~$s_{2}.$ There exist eight $\Z^2$-connected subsets of the
lattice~$\Z^2$ which satisfy the conditions of Proposition~\ref{suppthm},
namely
$$
\begin{array}{clcr}
S_1&=\{ (s_1,s_2)\in \Z^2: 1\leq s_1\leq 2,1\leq s_2\leq 3 \},\\
S_2&=\{ (s_1,s_2)\in \Z^2: 4\leq s_1, 5\leq s_2            \},\\
S_3&=\{ (s_1,s_2)\in \Z^2: 5\leq s_2, s_1+s_2 \leq 0       \},\\
S_4&=\{ (s_1,s_2)\in \Z^2: 4\leq s_1, s_1+s_2 \leq 0       \},\\
S_5&=\{ (s_1,s_2)\in \Z^2: 4\leq s_1, 1\leq s_2\leq 3      \},\\
S_6&=\{ (s_1,s_2)\in \Z^2: s_1+s_2 \leq 0, 1\leq s_2\leq 3 \},\\
S_7&=\{ (s_1,s_2)\in \Z^2: 1\leq s_1\leq 2, 5\leq s_2      \},\\
S_8&=\{ (s_1,s_2)\in \Z^2: 1\leq s_1\leq 2, s_1+s_2 \leq 0 \}.\\
\end{array}
$$
These irreducible supports of solutions to~(\ref{ex2diff}) are
displayed in Figure~1.
\vskip0.5cm
\hskip0.5cm
\begin{minipage}{7cm}
\begin{picture}(200,200)
   \put(100,10){\vector(0,1){180}}
   \put(10,100){\vector(1,0){180}}
   \put(180,88){\small $s_1$}
   \put(85,180){\small $s_2$}
   \put(130,20){\line(0,1){160}}
   \put(20,140){\line(1,0){160}}
   \put(20,180){\line(1,-1){160}}
   \put(20,130){\line(1,0){160}}
   \put(120,20){\line(0,1){160}}
   \put(110,110){\circle*{3}}
   \put(120,110){\circle*{3}}
   \put(110,120){\circle*{3}}
   \put(120,120){\circle*{3}}
   \put(110,130){\circle*{3}}
   \put(120,130){\circle*{3}}
   \put(140,150){\circle*{3}}
   \put(140,160){\circle*{3}}
   \put(140,170){\circle*{3}}
   \put(150,150){\circle*{3}}
   \put(150,160){\circle*{3}}
   \put(150,170){\circle*{3}}
   \put(160,150){\circle*{3}}
   \put(160,160){\circle*{3}}
   \put(160,170){\circle*{3}}
   \put(50,150){\circle*{3}}
   \put(40,150){\circle*{3}}
   \put(40,160){\circle*{3}}
   \put(30,150){\circle*{3}}
   \put(30,160){\circle*{3}}
   \put(30,170){\circle*{3}}
   \put(140,60){\circle*{3}}
   \put(140,50){\circle*{3}}
   \put(140,40){\circle*{3}}
   \put(150,50){\circle*{3}}
   \put(150,40){\circle*{3}}
   \put(160,40){\circle*{3}}
   \put(140,110){\circle*{3}}
   \put(140,120){\circle*{3}}
   \put(140,130){\circle*{3}}
   \put(150,110){\circle*{3}}
   \put(150,120){\circle*{3}}
   \put(150,130){\circle*{3}}
   \put(160,110){\circle*{3}}
   \put(160,120){\circle*{3}}
   \put(160,130){\circle*{3}}
   \put(90,110){\circle*{3}}
   \put(80,110){\circle*{3}}
   \put(80,120){\circle*{3}}
   \put(70,110){\circle*{3}}
   \put(70,120){\circle*{3}}
   \put(70,130){\circle*{3}}
   \put(60,110){\circle*{3}}
   \put(60,120){\circle*{3}}
   \put(60,130){\circle*{3}}
   \put(50,110){\circle*{3}}
   \put(50,120){\circle*{3}}
   \put(50,130){\circle*{3}}
   \put(40,110){\circle*{3}}
   \put(40,120){\circle*{3}}
   \put(40,130){\circle*{3}}
   \put(30,110){\circle*{3}}
   \put(30,120){\circle*{3}}
   \put(30,130){\circle*{3}}
   \put(110,150){\circle*{3}}
   \put(120,150){\circle*{3}}
   \put(110,160){\circle*{3}}
   \put(120,160){\circle*{3}}
   \put(110,170){\circle*{3}}
   \put(120,170){\circle*{3}}
   \put(110,90){\circle*{3}}
   \put(110,80){\circle*{3}}
   \put(120,80){\circle*{3}}
   \put(110,70){\circle*{3}}
   \put(120,70){\circle*{3}}
   \put(110,60){\circle*{3}}
   \put(120,60){\circle*{3}}
   \put(110,50){\circle*{3}}
   \put(120,50){\circle*{3}}
   \put(110,40){\circle*{3}}
   \put(120,40){\circle*{3}}
   \put(101,103){\small $S_1$}
   \put(170,175){\small $S_2$}
   \put(13,155){\small $S_3$}
   \put(145,23){\small $S_4$}
   \put(170,107){\small $S_5$}
   \put(13,108){\small $S_6$}
   \put(105,175){\small $S_7$}
   \put(105,23){\small $S_8$}
\end{picture}
\vskip-0.2cm

\noindent
{\scriptsize
{\bf Fig.~1} The irreducible supports of the solutions to the
Horn system~(\ref{ex2horn})
}
\end{minipage}
\hskip2cm
\begin{minipage}{3cm}
\vskip0.2cm
\begin{picture}(60,60)
   \put(0,0){\line(0,1){60}}
   \put(0,0){\line(1,0){60}}
   \put(60,0){\line(0,1){60}}
   \put(0,60){\line(1,0){60}}
   \put(40,0){\line(0,1){20}}
   \put(40,20){\line(-1,1){20}}
   \put(20,40){\line(-1,0){20}}
   \put(0,0){\circle*{3}}
   \put(20,0){\circle*{3}}
   \put(40,0){\circle*{3}}
   \put(0,20){\circle*{3}}
   \put(20,20){\circle*{3}}
   \put(40,20){\circle*{3}}
   \put(0,40){\circle*{3}}
   \put(20,40){\circle*{3}}
\end{picture}

\noindent
{\scriptsize
{\bf Fig.~2} The Newton po\-ly\-to\-pe of the resultant
of~(\ref{ex2horn})
}
\vskip0.4cm

\begin{picture}(60,60)
\put(30,30){\vector(0,1){30}}
\put(30,30){\vector(1,0){30}}
\put(30,30){\vector(0,-1){30}}
\put(30,30){\vector(-1,0){30}}
\put(30,30){\vector(1,1){30}}
\end{picture}

\noindent
{\scriptsize
{\bf Fig.~3} The fan of the Horn system~(\ref{ex2horn})
(see Section~\ref{fansec})
}

\end{minipage}
\vskip0.5cm
Using the formula~(12) in~\cite{Sadykov} for defining the periodic
function~$\phi(s),$ one can compute the sums of the corresponding
Laurent series. Let~$y_{i}(x)$ denote the series solution
to~(\ref{ex2horn}) with the support~$S_{i}.$ These functions are
defined up to unessential constant factors which we choose in a specific
way in order to make the formulas simpler. Computations (which were
performed in MAPLE) show that
$$
y_1(x)=3x_1x_2+4x_1x_{2}^{2}+3x_1x_{2}^{3}+3x_{1}^{2}x_2+6x_{1}^{2}x_{2}^{2}
+6x_{1}^{2}x_{2}^{3},
$$
$$
y_5(x)=x_{1}^4x_2(6x_{1}^{3}x_{2}^{2} +6x_{1}^{3}x_2
-27x_{1}^{2}x_{2}^{2} +3x_{1}^{3} -26x_{1}^{2}x_2  +
\phantom{--------}
$$
$$
\phantom{----}
45x_1x_{2}^2 -12x_{1}^2 + 40x_1x_2 - 30x_{2}^{2} +15x_1 -20x_2-6)/(1-x_1)^5,
$$
$$
y_7(x)=\frac{x_1x_{2}^5(6x_1x_{2}^2-18x_1x_2+3x_{2}^2+15x_1-8x_2+5)}
{(1-x_2)^4},
$$
$$
y_2(x)=\frac{x_1x_2(6x_{1}^2+14x_1x_2+5x_{2}^2 -9x_1 -8x_2
+3)}{(1-x_1-x_2)^4} - y_{1}(x) - y_{5}(x) + y_{7}(x)
$$
(we omit an explicit but cumbersome formula for~$y_{2}(x)$). The series
supported in $S_2,S_3,S_4$ represent the same solution to our system
since they represent the same rational function in different domains.
Finally, $y_{6}(x) = y_1(x)+y_5(x)$ and $y_{8}(x) = y_1(x)+y_7(x).$
It follows from Theorem~2.8 in~\cite{Sadykov} that the space of holomorphic
solutions to the system~(\ref{ex2horn}) has dimension~4 at any point
$x\in\C^2$ such that $(1-x_{1})(1-x_{2})(1-x_{1}-x_{2})\neq 0.$
Hence the rational functions $y_1(x),y_2(x),y_5(x),y_7(x)$ form a basis
in this space. Notice that the resultant of the principal symbols of the
operators in the system~(\ref{ex2horn}) is given by the polynomial
$(x_{1}x_{2})^4(1-x_{1})(1-x_{2})(1-x_{1}-x_{2}).$
\label{example1}
\end{example}

Recall that a convex cone is called {\it strongly convex} if it
does not contain any lines through the origin.
To conclude this section we formulate one more statement on the
properties of supports of hypergeometric series which will be
used in the sequel.

\begin{proposition}
A nonconfluent hypergeometric series with the support~$S$ has a
non\-empty domain of convergence if and only if the convex hull
of~$S$ is a polyhedral set which is contained in a translation
of a strongly convex cone. The domain of convergence of the
series~(\ref{solseries}) is independent on the parameters
$c_1,\ldots,c_p$ in the formula~(\ref{hypcoeff}) (we disregard
exceptional values of these parameters for which~(\ref{solseries})
terminates or reduces to a linear combination of hypergeometric
series in fewer variables).
\label{convergenceprop}
\end{proposition}
The first conclusion of this proposition follows from
Proposition~\ref{polyhedral}, the lemma in~\S~4.1 of~\cite{GGR}
and the properties of hypergeometric series in one variable
(see Chapter~1 in~\cite{Srivastava}). The second conclusion of
the proposition follows from Theorem~1 in~\S~4.1 of~\cite{Srivastava}.

  Finally we remark that there exists a simple relation between the
domain of convergence of a nonconfluent hypergeometric series and
its support. This relation is described by the two-sided Abel lemma
which is proved in Section~\ref{meromorphicsec}.


\section{The Fan of the Horn System
\label{fansec}}
\noindent
By an {\it affine convex cone} we mean a set of the form $C+\xi,$
where~$C$ is a convex cone in~$\R^{n}$ with apex at the origin
and $\xi\in\R^n.$ Let $C_1 + \xi_1$ and $C_2 + \xi_2$ be affine
convex cones with~$C_1,C_2$ being convex cones and
$\xi_1,\xi_2 \in\R^n.$ We say that $C_1 + \xi_1$ is
{\it smaller} than $C_2 + \xi_2$ if $C_1\subset C_2.$
If $C_1=C_2$ then the corresponding affine cones are said to be equal.
For a convex set $B\subset\R^n$ its {\it recession cone}~$C_{B}$ is
defined to be $C_{B} = \{ s\in\R^n : u + \lambda s \in B, \, \forall u\in B,
\lambda\geq 0\}$ (see \cite{Ziegler}, Chapter~1). That is, the recession
cone of a convex set is the maximal element in the family of those cones
whose shifts are contained in this set.

For the reason of brevity the recession cone of the
convex hull of the support of a Puiseux series solution to the Horn
system will be referred to as {\it the cone of its support.}
It has nonempty interior if and only if the corresponding hypergeometric
series cannot be represented as a linear combination of hypergeometric
series in fewer variables which depend monomially on the original ones.
In Example~\ref{example1} the cone of the irreducible support~$S_2$ is
the positive quadrant, the cone of~$S_5$ is
$\{ (s_1,s_2) : s_1\geq 0, \, s_2=0 \},$ the cone of~$S_1$ is the origin.

Here and later we assume that the rank of the matrix with the rows
$A_{1}, \ldots ,A_{p}$ is~$n$ since otherwise the series with the
coefficient~(\ref{hypcoeff}) can be reduced to a hypergeometric series
in fewer variables.
Let $I=(i_{1}, \ldots ,i_{n}),$ $i_{j}\in \{1, \ldots ,p\}$ be a multi-index
such that the vectors $A_{i_1}, \ldots ,A_{i_n}$ are linearly independent.
Let~$\gamma_{I}$ be the solution of the system of linear equations
$\langle A_{i_j}, s \rangle - c_{i_j} = 0,$ $j=1, \ldots ,n$
and define the set~$K_{I}$ by
$K_{I} = \{ s\in\Z^n :
\langle A_{i_j},s+\gamma_{I}\rangle - c_{i_j} \leq 0,
\, j=1, \ldots ,n \}.$
Let~$\Z^n +\gamma$ denote the shift in~$\C^n$ of the lattice~$\Z^n$ with
respect to the vector~$\gamma.$

\begin{definition}
\rm
We say that the parameter $c=(c_1, \ldots ,c_p)\in\C^p$ is {\it generic,}
if for any multi-index~$I$ as above none of the hyperplanes
$\langle A_{j},s+\gamma_{I} \rangle - c_{j} =0,$
$j\not\in \{i_{1}, \ldots ,i_{n}\}$
meets the shifted lattice $\Z^n + \gamma_{I}.$
\label{genericdef}
\end{definition}

\begin{proposition}
If the vector $c=(c_{1}, \ldots ,c_{p})$ is generic then there exists a
one-to-one correspondence between the $n$-dimensional
cones of the supports of the
convergent series solutions to the Horn system of the form~(\ref{solseries})
and the multi-indices $I=(i_{1}, \ldots ,i_{n})$ such that the vectors
$A_{i_1}, \ldots ,A_{i_n}$ are linearly independent. The recession cone of
the convex hull of the support of any such series is strongly convex
and polyhedral.
\label{polyhedralold}
\end{proposition}

\begin{proof}
For a multi-index~$I$ as above consider the shifted Laurent series
\begin{equation}
y_{I}(x) = \sum_{s\in K_I} t^s
\prod\limits_{i=1}^{p} \Gamma (\langle A_{i},s + \gamma_I \rangle - c_{i})
x^{s+\gamma_I}.
\label{convergent}
\end{equation}
Since the parameter~$c$ is assumed to be generic, it follows from
Proposition~\ref{suppthm} that the coefficient of the
series~(\ref{convergent}) satisfies the equations~(\ref{difference})
everywhere on~$\Z^n,$ i.e., that~(\ref{convergent}) is at least a
formal solution to
the Horn system~(\ref{horn}). By Proposition~\ref{convergenceprop} the
series~(\ref{convergent}) has a nonempty domain of convergence since its
support is contained in a strongly convex (and simplicial) affine cone.
Thus with any multi-index~$I$ as above one can associate the
$n$-dimensional cone~$C_I$
of the support of the series~(\ref{convergent}).

Since we are interested in $n$-dimensional cones of the supports
of the series solutions to~(\ref{horn}), we do not consider
polynomial solutions to this system (which may exist even if the
parameters are generic).
It follows by Proposition~\ref{suppthm} that if the support of a formal
series solution to~(\ref{horn}) meets at most $n-1$ linearly independent
hyperplanes of the form $\langle A_{j}, s + \gamma \rangle - c_{j} = 0$
for some $\gamma\in\C^{n}$ then it cannot
be contained in any strongly convex affine cone and by
Proposition~\ref{convergenceprop} the series is divergent. By the
assumption the parameter~$c$ is generic and hence the support of
such a series cannot meet more than~$n$ hyperplanes of this form.
If it meets exactly~$n$ hyperplanes with the linearly independent
normals $A_{i_1}, \ldots ,A_{i_n}$ then the cone of the support of
this series must coincide with~$C_I$ since it is bounded by the same
hyperplanes. Thus the correspondence between linearly independent
subsets of the set of vectors $\{ A_{1}, \ldots ,A_{p} \}$ and the
$n$-dimensional cones
of the supports of shifted Laurent series solutions to~(\ref{horn})
is one-to-one. The claim about the recession cone of the convex hull
of the support of~$y_{I}(x)$ follows from Proposition~1.12
in~\cite{Ziegler} since the convex hull of~$K_{I}$ is a strongly
convex affine polyhedral cone.~\hfill~$\square$
\end{proof}

\begin{remark}
\rm
Proposition~\ref{polyhedralold} shows that adding new elements to the
family of vectors ${\{A_i\}}_{i=1}^{p}$ can only increase the number
of series solutions to the Horn system which is defined by the Ore-Sato
coefficient~(\ref{hypcoeff}) as long as the vector~$c$ remains generic.
\label{refinerem}
\end{remark}

We now associate with a nonconfluent Horn system a set of strongly convex
polyhedral cones which will play an important role in the sequel.
Recall that for a cone $C\subset\R^n$ its dual is defined by
$C^{\vee}= \{ v\in\R^n : \langle u,v \rangle \geq 0,\, \forall u\in C\}.$
For any multi-index $I=(i_{1}, \ldots ,i_{n})$ such that the vectors
$A_{i_1}, \ldots ,A_{i_n}$ are linearly independent we denote by~$C_{I}$
the recession cone of the convex hull of the set~$K_I$ whose shift
supports the series~(\ref{convergent}). We partially order the finite
family~$\{ C_{I} \}$ of strongly convex polyhedral cones with respect
to inclusion and denote the maximal elements by
$C_{I^{(1)}}, \ldots ,C_{I^{(d)}}.$ Let us introduce the cones
$B_j= -C_{I^{(j)}}^{\vee},$ $j=1, \ldots ,d.$ Since for any~$I$ as above the
polyhedral cone~$C_{I}$ has a nonempty interior, it follows that~$B_j$
is a strongly convex polyhedral cone. The nonconfluency
condition~(\ref{nonconfluency}) implies that
$\bigcup_{j=1}^{d} B_j = \R^n.$ If the cones $B_1, \ldots ,B_d$ can be
identified with the set of the maximal cones of some complete fan then
we call it {\it the fan of the Horn system~(\ref{horn}).}

If $n=2$ then ${\{ B_j \}}_{j=1}^d$ is always the set of the maximal
cones of some complete fan. For $n\geq 3$ this is not necessarily the
case. For instance, let $n=3$ and let $A_1=(1,0,0),$ $A_2=(0,1,0),$
$A_3=(0,0,2),$ $A_4=(-1,0,-1),$ $A_5=(0,-1,-1).$ The multi-indices
$I^{(1)}=(1,4,5)$ and $I^{(2)}=(2,4,5)$ define maximal cones but the
intersection of their duals has a nonempty interior.


\section{Minimality of the singularities of \newline
hypergeometric functions and discriminants
\label{amoebasec}}

As we have already mentioned in the introduction, the singular set of
a hypergeometric function in one variable is minimal in th esense that
its amoeba consists of a single point. In this section
we will prove that multivariate rational hypergeometric functions enjoy
an analogous property. It turns out to be convenient to express this
property using the notion of amoebas which was introduced by Gelfand
et al. in~\cite{GKZ} (see Chapter~6,~\S~1).
The {\it amoeba}~$\mathcal{A}_f$ of a Laurent polynomial~$f(x)$
(or of the algebraic hypersurface $f(x)=0$) is defined to be the image of
the hypersurface~$f^{-1}(0)$ under the map
${\rm Log } : (x_1,\ldots,x_n)\mapsto (\log |x_1|,\ldots,\log |x_n|).$
This name
is motivated by the typical shape of~$\mathcal{A}_f$ with tentacle-like
asymptotes going off to infinity (see Figure~5). We quote the following
result from~\cite{GKZ} (see Chapter~6, Corollary~1.6), which describes
the connection between the amoeba of a Laurent polynomial~$f$ and Laurent
series developments of~$1/f.$

\vskip0.1cm
\noindent
{\bf Theorem A}
{\rm (Gelfand, Kapranov, Zelevinsky~\cite{GKZ})}
{\it
The connected components of the amoe\-ba complement~$^c\!\mathcal{A}_{f}$
are convex, and they are in bijective correspondence with the different
Laurent series expansions centered at the origin of the rational
function~$1/f.$
} 
\vskip0.1cm

Recall that the Newton polytope~$\mathcal{N}_f$ of a Laurent
polynomial~$f$ is defined to be the convex hull in~$\R^n$ of the support
of~$f.$ The following result shows that the Newton polytope~$\mathcal{N}_f$
reflects the structure of the amoeba~$\mathcal{A}_f$ (see Theorem~2.8
and Proposition~2.6 in~\cite{FPT}).

\vskip0.1cm
\noindent
{\bf Theorem B}
{\rm (Forsberg, Passare, Tsikh \cite{FPT})}
{\it
Let~$f$ be a Laurent polynomial and let~$\{M\}$ denote the family
of connected components of the amoeba complement~$^c\!\mathcal{A}_f.$
There exists an injective function
$\nu : \{M\}\rightarrow \Z^n \cap \mathcal{N}_f$ such  that
the cone which is dual to~$\mathcal{N}_f$ at the point~$\nu(M)$
coincides with the recession cone of~$M.$
} 
\vskip0.1cm

The cited theorems imply that the number of Laurent series expansions
with the center at the origin of the rational function~$1/f$ is at least
equal to the number of vertices of the Newton polytope~$\mathcal{N}_f$
and at most equal to the number of integer points in~$\mathcal{N}_f.$
Varying the coefficients of the Laurent polynomial~$f$ with the fixed
Newton polytope~$\mathcal{N}_f,$ one can attain the upper
(see~\cite{Mikhalkin}) as well as the lower (see~\cite{Rullgard})
bounds for the number of connected components of~$^c\!\mathcal{A}_f.$
Moreover, the vertices of the Newton polytope are always assumed
by the function~$\nu$ and by Theorem~B the recession
cones of those connected components of~$^c\!\mathcal{A}_f$ which
correspond to the vertices of~$\mathcal{N}_f$ have nonempty interior.

In this section we show that if~$f$ is the defining polynomial of
the singular locus of a hypergeometric function
then the number of connected components of~$^c\!\mathcal{A}_f$
equals the number of vertices of~$\mathcal{N}_f.$
For the sake of brevity we use the following definition.

\begin{definition}
The amoeba~$\mathcal{A}_{f}$ of a Laurent polynomial~$f$ (or,
equivalently, the algebraic hypersurface $f(x)=0$) is called
{\rm solid} if the number of connected components of the amoeba
complement~$^c\!\mathcal{A}_{f}$ equals the number of vertices of the
Newton polytope~$\mathcal{N}_{f}.$
\label{minimaldef}
\end{definition}
In view of Theorem~B it is obvious that the amoeba~$\mathcal{A}_{f}$
is solid if and only if the recession cone of every connected
component of the set~$^c\!\mathcal{A}_{f}$ has a nonempty interior.
The main observation in this section is the following theorem.
\begin{theorem}
The singular hypersurface of any nonconfluent hypergeometric function
has a solid amoeba.
\label{rathypmin}
\end{theorem}
\begin{proof}
Let~$\mathcal{A}$ be the amoeba of the resultant of the Horn system
(as defined in Section~\ref{resultantsec}) and let
$M\subset\,  ^c\!\mathcal{A}$ be a connected component of its complement.
By the remark after Definition~\ref{minimaldef} it suffices to show that
the recession cone~$C_{M}$ of the set~$M$ has a nonempty interior.

Recall that in this paper we only deal with Horn systems satisfying
the assumptions made in Section~\ref{resultantsec}.
The condition that the projection of the characteristic variety of
the Horn system onto the variable space is its proper algebraic
subset implies that the Horn system in question is holonomic
(see Chapter~3 of~\cite{Bjork1}).
Hence it has finitely many analytic solutions in a neighbourhood of its
nonsingular point.

Our next argument was inspired by the proof of Theorem~2.4.12 in~\cite{SST}.
Let $y_{1},\ldots,y_{r}$ be a basis in the space of holomorphic
solutions to~(\ref{horn}) on a simply connected domain
in~${\rm Log}^{-1} M.$ Recall that~$J$ denotes the ideal generated
by the differential operators in the Horn system.
Let $\{ 1, \partial^{\alpha(1)},\ldots, \partial^{\alpha(r-1)} \}$
be a basis of the quotient
$\C(x) \langle \partial \rangle / \C(x) \langle \partial \rangle J,$
where $\partial =$ $(\partial_{1}, \ldots, \partial_{n})=$
$\frac{\partial}{\partial x_{1}},\ldots, \frac{\partial}{\partial x_{n}}$
and $\C(x) \langle \partial \rangle = $
$\C(x_{1},\ldots,x_{n}) \langle \partial_{1},\ldots, \partial_{n} \rangle$
is the algebra generated by polynomials in
$\partial_{1},\ldots,\partial_{n}$ and rational functions in
$x_{1},\ldots,x_{n}.$
Put
$$
\Phi(x) =
\left(
\begin{array}{ccc}
y_{1}                       & \ldots &   y_{r}                       \\
\partial^{\alpha(1)}y_{1}   & \ldots &   \partial^{\alpha(1)}y_{r}   \\
\ldots                      & \ldots &   \ldots                      \\
\partial^{\alpha(r-1)}y_{1} & \ldots &   \partial^{\alpha(r-1)}y_{r} \\
\end{array}
\right).
$$
Since~$\{ y_{i} \}$ is a basis, it follows that
${\rm det} (\Phi)\not\equiv 0$ and~$\Phi$ is a (matrix-valued)
multi-valued holomorphic function on~${\rm Log}^{-1} M.$
By Theorem~A the set~$M$ is convex and hence~${\rm Log}^{-1} M$ is
a Reinhardt domain with the center at the origin. Its fundamental
group~$\pi_{1}({\rm Log}^{-1} M)$ is isomorphic to the direct product
of the fundamental groups of at most~$n$ punched disks with the
center at the origin. Thus~$\pi_{1}({\rm Log}^{-1} M)$ is a free
Abelian group generated by the elements~$\eta_{i}$ which
encircle~$x_{i}=0$ (some of these elements might be trivial).

Consider the analytic continuation~$\eta_{i}^{*}\Phi$ of the matrix~$\Phi$
along the path~$\eta_{i}.$
Since the first row of~$\eta_{i}^{*}\Phi$ is again a basis of solutions,
there exists an invertible matrix~$V_{i},$ which is called the monofromy
matrix, satisfying $\eta_{i}^{*} \Phi = \Phi V_{i}.$
Since~$\pi_{1}({\rm Log}^{-1} M)$ is Abelian, the matrices~$V_{i}$
commute with one another. Hence there exists a commutative family
of matrices~$W_{i}$ such that $e^{2 \pi \sqrt{-1}\, W_{i}} = V_{i}.$
Define the matrix
$$
\Psi (x):= \Phi(x) x_{1}^{-W_{1}} \ldots x_{n}^{W_{n}}.
$$
The monodromy of~$\Phi(x)$ is killed
by~$x_{1}^{-W_{1}} \ldots x_{m}^{-W_{m}}$ since
$\eta_{i}^{*} x_{i}^{-W_{i}} = V_{i}^{-1} x_{i}^{-W_{i}}.$
Hence~$\Psi(x)$ is a single-valued function on~${\rm Log}^{-1} M.$
By Lemma~2 in Chapter~4 of~\cite{Bolibrukh} any solution to the Horn
system in the domain~${\rm Log}^{-1} M$ can be written as a
polynomial in Puiseux monomials and~${\rm log }\, x_{i}$ with single-valued
coefficients. Here by a Puiseux monomial we mean a monomial with
arbitrary (complex) exponent vector.

Let us write such a solution in the form
$y(x) = \sum\limits_{\alpha,\beta} h_{\alpha \beta}(x) x^{\alpha}
({\rm log }\, x)^{\beta},$ where~$h_{\alpha \beta}(x)$ are single-valued
functions in~${\rm Log}^{-1} M,$
$({\rm log }\, x)^{\beta} := ({\rm log }\, x_{1})^{\beta_{1}} \ldots
({\rm log }\, x_{n})^{\beta_{n}}$ and the sum is finite.
Let~$\beta_{1}^{'}$ be the highest power of~${\rm log }\, x_{1}$
appearing in the expression for~$y(x).$
Any single-valued function in a Reinhardt domain can be expanded into
a Laurent series.
Expanding the functions~$h_{\alpha \beta}$ into Laurent series and
computing the action of the operators in the Horn system on~$y(x),$
we conclude that the coefficients of the expansion for~$h_{\alpha \beta}$
satisfy difference relations of the form~(\ref{difference}).
The first of these relations yields an ordinary hypergeometric
differential equation for the restriction of~$y(x)$ to a suitable line.
It is known that no logarithms may appear in a solution to an
ordinary generalized hypergeometric differential equation with generic
parameters (see~\cite{Evgrafov}). By induction over the highest power
of~${\rm log }\, x_{1}$ appearing in the expression for~$y(x)$ we
conclude that~${\rm log }\, x_{1}$ does not appear at all if the
parameters of the Horn system are sufficiently general.
By the symmetry of the variables it follows that any solution to
a Horn system with generic parameters in the domain~${\rm Log}^{-1} M$
can be represented as a Puiseux series.

For $\zeta\in\partial M$ let $Y_{\zeta}\subset\R^n$ denote the
half-space which is bounded by a supporting hyperplane of~$M$ at
the point~$\zeta$ and contains~$M.$ There exists a sequence of
points ${\{\zeta_i\}}_{i=1}^{\infty} \subset \partial M$ such that
the recession cone of the set
$\bigcap_{i=1}^{\infty} Y_{\zeta_i}$ coincides with~$C_{M}.$
Since~$\mathcal{A}$ is the logarithmic image of the set of singularities
of the function~$y(x),$ for any $i\in\N$ there exists a
germ~$\mathcal{G}_i$ of~$y(x)$ which cannot be continued analytically
through at least one point in the fiber ${\rm Log }^{-1}\zeta_i.$
As we have remarked earlier, the analytic
continuation of~$\mathcal{G}_i$ into the domain ${\rm Log}^{-1} M$
can be expanded into a Puiseux series~$L_{i}$ whose domain of
convergence contains ${\rm Log}^{-1} M.$ Let
$L^{(k)}=\sum_{i=1}^{k} L_{i}.$ The series~$L^{(k)}$
satisfies the same hypergeometric system of equations as~$y(x)$ since
it is a linear combination of solutions to this system. By the
construction~$L^{(k)}$ is not identically equal to zero.
We denote the domain of convergence of the
series~$L^{(k)}$ by~$\Omega_{k}.$
By the construction $M\subset {\rm Log }\,\Omega_{k}$ and hence~$\Omega_{k}$
is nonempty. Moreover the recession cone $C_{{\rm Log }\Omega_{k}}$
is a subset of the recession cone of the finite intersection
$\bigcap_{i=1}^{k} Y_{\zeta_i}.$

Suppose that the cone~$C_{M}$ has the empty interior. The two-sided Abel
lemma which will be proved in Section~\ref{meromorphicsec} states that
for a nonconfluent hypergeometric Puiseux series~$L$ with the domain of
convergence~$\Omega$ one has $C_{{\rm Log }\, \Omega} = - C_{L}^{\vee},$
where~$C_{L}$ is the cone of the support of~$L$ and
$C_{{\rm Log }\, \Omega}$ is the recession cone of the set
${\rm Log }\, \Omega.$ Thus we have
$-C_{L^{(k)}}^{\vee} = C_{{\rm Log }\Omega_{k}} \subset
C_{\bigcap_{i=1}^{k} Y_{\zeta_i}}$
and hence the set
$\bigcup_{k=1}^{\infty} C_{L^{(k)}}$ is not strongly convex.
By Proposition~\ref{polyhedral} the cone~$C_{L^{(k)}}$ is polyhedral
with its boundary being a subset of the union of the zero sets of the
polynomials $P_1,\ldots,P_n,$ $Q_1,\ldots,Q_n.$ Since this union is a finite
arrangement of hyperplanes it follows that the family of
cones~${\{C_{L^{(k)}}\}}_{k=1}^{\infty}$
can only contain a finite number of distinct elements. Therefore
there exists $m\in\N$ such that the cone~$C_{L^{(m)}}$ is not strongly
convex. This contradicts the statement of Proposition~\ref{convergenceprop}
and completes the proof.~\hfill~$\square$
\end{proof}

Let us recall the definition of $\mathcal{A}$-discriminant
which was introduced by Gelfand, Kapranov and Zelevinsky
(see~\cite{Kapranov}).
Let~$\mathcal{A}$ be a finite subset of~$\Z^n$ and let~$f$
be a generic polynomial with the support~$\mathcal{A},$ i.e.,
$f=\sum\limits_{\alpha\in\mathcal{A}} c_{\alpha}x^{\alpha}.$
The corresponding $\mathcal{A}$-discriminant is defined to be
the polynomial in the coefficients~$c_{\alpha}$ which vanishes
whenever~$f$ together with all of its partial derivatives have a common
zero.

A hypergeometric function satisfying the Gelfand-Kapranov-Zelevinsky
system of equations (see~\cite{GKZ0}) has singularities along the
zero locus of the corresponding $\mathcal{A}$-discriminant.
There always exists a monomial change of variables which transforms an
$\mathcal{A}$-hypergeometric series into a Horn series (see Section~2
in~\cite{Kapranov}).
This monomial change of variables corresponds to a linear
transformation of the amoeba space and hence it cannot affect
the solidness of an amoeba.
(More precisely, the preimage of any point in the amoeba space under
this mapping is an affine subspace and hence the preimage of
a solid amoeba is also solid.)
Using Theorem~\ref{rathypmin} we arrive at the following corollary.

\begin{corollary}
The zero set of any $\mathcal{A}$-discriminant has a solid amoeba.
\label{adiscmin}
\end{corollary}

Theorem~\ref{rathypmin} allows also to derive the following property of
the classical discriminant of the general algebraic equation
$y^m + c_{1}y^{m_1} + \cdots + c_{n}y^{m_n} + c_{n+1}=0,$
where $m,m_i\in\N,$ $m > m_1 > \ldots > m_n \geq 1,$
$y$~is the unknown. We provide the following corollary with a proof
since the solution to a general algebraic equation satisfies a
system of differential equations which is slightly different
from~(\ref{horn}).

\begin{corollary}
The amoeba of the discriminant of a general algebraic equation is solid.
\label{discriminant}
\end{corollary}
\begin{proof}
By a monomial change of the variable~$y$ and the coefficients
$c_{1}, \ldots, c_{n+1}$
any algebraic equation can be reduced to an equation of the form
\begin{equation}
y^{m} + x_{1}y^{m_1} + \ldots + x_{n}y^{m_n} - 1 = 0,
\label{algeq}
\end{equation}
where $x=(x_1,\ldots,x_n)\in\C^n.$
It was shown in~\cite{Mellin} that the solution~$y(x)$ to~(\ref{algeq})
(which is considered as a multi-valued analytic function depending on
$x_1,\ldots,x_n$) satisfies the system of partial differential equations
$$
{(-1)}^{m_i}m^{m} \frac{\partial^{m}y}{\partial x_{i}^{m}} =
\prod_{j=0}^{m_i -1}
(m_1\theta_1 + \ldots + m_n\theta_n + 1 + mj) \times
\phantom{----------}
$$
\begin{equation}
\phantom{---}
\times
\prod_{j=0}^{m_{i}^{'}-1}
(m_{1}^{'}\theta_1 + \ldots + m_{n}^{'}\theta_n - 1 + mj)y,
\,\,\, i=1,\ldots,n,
\label{mellinsyst}
\end{equation}
where $m_{i}^{'}=m-m_{i}.$
The singular set of the function~$y(x)$ is the discriminant of the
equation~(\ref{algeq}).
Multiplying the $i$th equation of~(\ref{mellinsyst}) with~$x_{i}^{m},$
using the identity
$x_{i}^{m}\frac{\partial^{m}}{\partial x_{i}^{m}}=
\prod_{j=0}^{m-1}(\theta_i - j)$ and making the monomial change of
variables $\xi_i = x_{i}^{m}$ we reduce the Mellin system~(\ref{mellinsyst})
to a system of the form~(\ref{horn}). Thus~$y(\xi)$ is a nonconfluent
hypergeometric
function in the sense of Horn. Since the function~$y(x)$ is of finite
branching, so is~$y(\xi).$ By Theorem~\ref{rathypmin} the singular set
of~$y(\xi)$ has a solid amoeba. Since a monomial change of variables 
corresponds
to a linear transformation of the amoeba space (see~\cite{FPT}), it follows
that such a change of variables cannot affect the solidness of the
singularity of~$y(x).$ Thus the amoeba of the discriminant of the algebraic
equation~(\ref{algeq}) is solid.~\hfill~$\square$
\end{proof}

The cubic equation is considered in detail in Example~\ref{Cardano}.
The amoeba of the singular locus of a solution to the reduced system
is displayed in Figure~5.

Theorem~\ref{rathypmin} implies in particular that the number of
connected components of the complement of the amoeba of the singular
hypersurface of a rational hypergeometric function equals the number
of vertices of the Newton polytope of its denominator. It turns out
that in some cases knowing the hypergeometric system which is
satisfied by a given rational function allows one to compute the
number of vertices of the Newton polytope of its denominator.
We illustrate this fact by means of the following important family
of rational hypergeometric functions which are defined as the Bergman
kernels of complex ellipsoidal domains (see~\cite{FH} and~\cite{Zinov'ev}).
This family will be used in Section~\ref{classificationsec} for
describing rational hypergeometric functions satisfying some systems
of equations of the Horn type.

Consider the family of {\it complex ellipsoidal domains} defined by
$$
D^{p_{1},\ldots,p_{n}}=
\{ x\in\C^n : {|x_{1}|}^{2/p_{1}}+\ldots+{|x_{n}|}^{2/p_{n}} < 1 \},
$$
where $p_{i}=1,2,3,\ldots,$ $i=1,\ldots,n.$ The Bergman kernel
$K_{p_{1},\ldots,p_{n}}(x)$ for this domain was computed
explicitly in~\cite{Zinov'ev}. It was shown that
\begin{equation}
K_{p_{1},\ldots,p_{n}}(x)=
\frac{1}{\pi^{n}} \sum\limits_{s\in\N_{0}^n}
\frac{\Gamma(p_{1}(s_{1}+1)+ \ldots +p_{n}(s_{n}+1) + 1)}
{\prod_{i=1}^n p_{i}\Gamma(p_{i}(s_{i}+1))}x^s.
\label{bergman}
\end{equation}
The sum of this series is given by the function
\begin{equation}
K_{p_{1},\ldots,p_{n}}(x) = \frac{1}{\pi^n}\frac{1}{p_{1}\ldots p_{n}}
\frac{\partial^n}{\partial x_{1}\ldots \partial x_{n}}
\sum_{j_{1}=1}^{p_{1}}\ldots \sum_{j_{n}=1}^{p_{n}}
\frac{1}{1-y_{j_{1}1}-\ldots- y_{j_{n}n}},
\label{rational}
\end{equation}
where $y_{j_{i}i}=x_{i}^{1/p_{i}}\varepsilon_{j_{i}i},$\,
$\varepsilon_{j_{i}i}$ are all the $p_{i}$-th roots of unity,
$j_{i}=1,\ldots,p_{i},$ $i=1,\ldots,n.$ The expression under the sign
of the partial derivatives in~(\ref{rational}) was proved
in~\cite{Zinov'ev} to be rational in $x_{1},\ldots,x_{n}$ and to
have integral coefficients for any choice of $p_{1},\ldots,p_{n}$.
Let~$f_{p_{1},\ldots,p_{n}}$ denote the denominator of the rational
function~(\ref{rational}) (we normalize the denominator so that the
greatest common divisor of its coefficients equals~$1$). Our aim is
to find the number of connected components of the amoeba
complement~$^c\!\mathcal{A}_{f_{p_{1},\ldots,p_{n}}}.$
For any fixed vector $\gamma\in \C^n,$
${\rm Re}\,\gamma_{i} \in [0,1),$ there exist finitely many subsets
of the shifted lattice~$\Z^n +\gamma$ which satisfy the conditions in
Proposition~\ref{suppthm} and are contained in some strongly convex
affine cone. We call them {\it $\gamma$-admissible sets associated
with~(\ref{horn})}. A set is said to be {\it admissible} if it is
$\gamma$-admissible for some~$\gamma.$
\begin{proposition}
The number of connected components of the amoeba complement
$^c\!\mathcal{A}_{f_{p_{1},\ldots,p_{n}}}$ of the denominator
of the Bergman kernel~$K_{p_{1},\ldots,p_{n}}(x)$ equals~$n+1.$
\label{minimalprop}
\end{proposition}
\begin{remark}
\rm
The conclusion of Proposition~\ref{minimalprop} can be deduced from
Proposition~4.2 in~\cite{FPT} in the following way. Let us introduce
new variables $\xi_{i}=x_{i}^{1/p_{i}}.$
It follows from Proposition~4.2 in~\cite{FPT} that for any choice of
the indices $j_1\in\{1,\ldots,p_1\},\ldots,j_n\in\{1,\ldots,p_n\}$
the amoeba of the first-order polynomial
$1-\varepsilon_{j_{1}1}\xi_{1} - \ldots - \varepsilon_{j_{n}n}\xi_{n}$
is the same. By Corollary~4.5 in~\cite{FPT} the number of connected
components of its complement equals $n+1.$ Since a monomial change of
the variables $x_1,\ldots,x_n$ corresponds to a linear transformation of
the amoeba space (see~\cite{FPT}), it follows that the number of the
connected components of the complement of the amoeba of $f_{p_1,\ldots,p_n}$
also equals $n+1.$ This shows in particular that the amoeba of
$f_{p_1,\ldots,p_n}$ is solid.
\label{anotherproof}
\end{remark}

We give here another proof of Proposition~\ref{minimalprop} which only
uses hypergeometric properties of the Bergman kernels and does not use
the explicit form of their denominators.

\noindent {\sc Proof of Proposition~\ref{minimalprop}.}
The Newton polytope of~$f_{p_{1},\ldots,p_{n}}$ has nonzero
$n$-dimensional volume. Indeed, the restriction of
$K_{p_{1},\ldots,p_{n}}(x)$ to the complex line
$x_{1}=\ldots [i] \ldots =x_{n}=0$ is a rational function whose denominator
is given by $(1-x_{i})^{k_{i}},\,\,$ $k_{i}>0,$ $i=1,\ldots,n.$
(Here~$[i]$ is the sign of omission.) It follows by
Theorem~B that the number of connected components of
the amoeba complement~$^c\!\mathcal{A}_{f_{p_{1},\ldots,p_{n}}}$ cannot
be smaller than~$n+1.$

Let~$\varphi(s)$ denote the coefficient of the series~(\ref{bergman}),~i.e.,
$$
\varphi(s)=\frac{\Gamma(p_{1}(s_{1}+1)+\ldots+p_{n}(s_{n}+1)+1)}{
\prod\limits_{i=1}^n p_{i} \Gamma(p_{i}(s_{i}+1))}.
$$
Since for any $i=1,\ldots,n$ the function~$\varphi(s)$ satisfies the
equation
$$
\varphi(s+e_{i})\prod_{j=0}^{p_{i}-1}(p_{i}(s_{i}+1)+j)=
\varphi(s)\prod_{j=1}^{p_{i}}(p_{1}(s_{1}+1)+\ldots+p_{n}(s_{n}+1)+j),
$$
it follows that~$K_{p_{1},\ldots,p_{n}}(x)$ is a solution to the
following system of the Horn type
$$
x_{i}\left(
\prod_{j=1}^{p_{i}}
(p_{1}(\theta_{1} + 1) +\ldots+ p_{n}(\theta_{n} + 1 ) + j)
\right) K_{p_{1},\ldots,p_{n}}(x)=
\phantom{----------}
$$
\begin{equation}
\phantom{------}
\left(
\prod_{j=0}^{p_{i}-1}
(p_{i}\theta_{i} + j)
\right) K_{p_{1},\ldots,p_{n}}(x), \,\,\, i=1,\ldots,n.
\label{bergmanhorn}
\end{equation}
The number of irreducible $0$-admissible sets associated with the
system~(\ref{bergmanhorn}) equals $n+1.$ These sets are
$S_{0}=\N_{0}^n$ and $S_{i}=\{ s\in\Z^n :
p_{1}(s_{1}+1)+\ldots+p_{n}(s_{n}+1)+1 \leq 0,\,\, s_{j}\geq 0, \,\,
j\neq i \},$ $i=1,\ldots,n.$ (Notice that~(\ref{bergman}) is supported
in the $0$-admissible set~$\N_{0}^n.$) Since any expansion of a rational
solution to a Horn system into a Laurent series with the center at the
origin is supported in an irreducible $0$-admissible set, it follows that
the number of connected components of the amoeba
complement~$^c\!\mathcal{A}_{f_{p_{1},\ldots,p_{n}}}$ cannot exceed~$n+1.$
We have proved earlier that the Newton polytope of~$f_{p_{1},\ldots,p_{n}}$
has at least~$n+1$ vertices. Thus it follows from
Theorem~B that the number of connected components
of~$^c\!\mathcal{A}_{f_{p_{1},\ldots,p_{n}}}$ cannot be smaller than~$n+1$
and hence equals~$n+1.$ The proof is complete.~\hfill~$\square$
\vskip0.15cm
\begin{example}
\rm
Let $n=2,$ $p_{1}=3,$ $p_{2}=2.$ The denominator of the Bergman kernel
of the domain~$D^{3,2}$ is given by
$$
f_{3,2}(x)={(1 - 2x_{1} - 3x_{2} + x_{1}^2  - 6x_{1}x_{2} + 3x_{2}^2 -
x_{2}^3)}^3.
$$
By Proposition~\ref{minimalprop} the number of connected components
of the amoeba complement~$^c\!\mathcal{A}_{f_{3,2}}$ equals~3.
\label{exbergman}
\end{example}

The Bergman kernel~(\ref{bergman}) gives an example of a rational
hypergeometric function. The problem of describing the class of rational
hypergeometric functions was studied in~\cite{CDD} and~\cite{CDS}.
Observe however, that the definition of a hypergeometric function used
in these papers is based on the Gelfand-Kapranov-Zelevinsky system of
differential equations~\cite{GGR} rather than the Horn system.


\section{Meromorphic nonconfluent hypergeometric \newline
functions are rational
\label{meromorphicsec}}

The aim of this section is to show that a nonconfluent Horn
system~(\ref{horn}) cannot possess a meromorphic solution different
from a rational function (Theorem~\ref{merorat}).

The relation between the support of a general Puiseux series and its
domain of convergence is described by the Abel lemma (see Introduction
and~\cite{GKZ0}, \S~1). For hypergeometric series the following stronger
version of this statement holds.
\begin{lemma}
{\rm (Two-sided Abel's lemma)}
Suppose that a nonconfluent hypergeometric Puiseux series with
the support~$S$ has nonempty domain of convergence~$D.$ Let~$C$ be
the cone of~$S.$ Then for any $x^{(0)}\in D$ and for some
$x^{(1)}\in\C^n \setminus D$
$$
{\rm Log \,}(x^{(0)}) - C^{\vee} \subset {\rm Log \,}(D) \subset
{\rm Log \,}(x^{(1)}) - C^{\vee}.
$$
\label{twosidedAbel}
\end{lemma}
\begin{proof}
Let $y(x)=\sum_{s\in S}\varphi(s)x^s$ be a nonconfluent
hypergeometric Puiseux series.
The first inclusion follows from the general Abel lemma (see
Introduction). Let us prove the second inclusion.
Let $M\subset\R^n$ be the lattice generated by the elements of
the set~$S.$
By Proposition~\ref{convergenceprop} the domain~$D$ is independent on
the parameters $c_{1},\ldots,c_{p}$ of the coefficient~(\ref{hypcoeff})
as long as they remain generic. Thus we may without loss of
generality assume that $S=C\cap M.$ Since~$D$ is nonempty, it follows
by Proposition~\ref{convergenceprop} that~$C$ is a strongly convex
polyhedral cone.
Let $u^{(1)},\ldots,u^{(N)}\in M$ denote the generators of~$C,$ i.e.,
$C=\{ \lambda_{1}u^{(1)}+\ldots+\lambda_{N}u^{(N)} :
\lambda_{j}\geq 0, \, j=1,\ldots,N \}.$  For each $j=1,\ldots,N$ we
consider the restricted series $y_{j}(x) = \sum_{k=0}^{\infty}
\varphi(ku^{(j)}) x^{ku^{(j)}}.$
The nonconfluency condition~(\ref{nonconfluency}) implies that
$\sum_{i=1}^{p} \langle A_{i},u^{(j)} \rangle = 0.$
By the result on convergence of the generalized hypergeometric series
in one variable (see~\cite{GGR}, \S~1.1) the domain of convergence
of~$y_{j}(x)$ is contained in the set
$\{ x\in\C^n : |x^{u^{(j)}}| < r_{j}\}$ for some constant $r_{j}>0.$
This shows that
${\rm Log \,}(D) \subset \{ v\in\R^n : \langle u^{(j)}, v \rangle <
\log r_{j}, \, j=1,\ldots,N \}.$
Since~$C$ is strongly convex, we can choose $\xi\in\R^n$ such that
$m_{j}:= \langle u^{(j)},\xi \rangle >0.$ Let
$$
x^{(1)} \in {\rm Log}^{-1}
\left( \xi\max_{j=1,\ldots,N} \frac{\log r_{j}}{m_{j}}\right),
$$
then $\langle u^{(j)}, {\rm Log \,}x^{(1)} \rangle \geq \log r_{j},$
$j=1,\ldots,N$
and hence ${\rm Log \,}(D)\subset \{ v\in\R^n : \langle u^{(j)},
v-{\rm Log \,} x^{(1)} \rangle \leq 0, \, j=1,\ldots,N\} =$
${\rm Log \,} x^{(1)} - C^{\vee}.$ The proof is complete.~\hfill~$\square$
\end{proof}

The two-sided Abel lemma enables us to prove the following theorem which
is the main result in this section.
\begin{theorem}
Any meromorphic nonconfluent hypergeometric function is rational.
\label{merorat}
\end{theorem}
\begin{proof}
Let~$y(x)$ be a meromorphic nonconfluent hypergeometric function.
By definition~$y(x)$ is a solution to the Horn system~(\ref{horn}).
Since~$y(x)$ is nonconfluent, it follows by
Proposition~\ref{convergenceprop} and the two-sided Abel lemma
(Lemma~\ref{twosidedAbel}) that the domain of convergence of any
shifted Laurent series representing~$y(x)$ is not all of~$(\C^{*})^n.$
Therefore, using the assumption that~$y(x)$ is meromorphic, we can write
it in the form~$h(x)/g(x),$ where~$h(x)$ is entire and~$g(x)$ is some
polynomial which is not a monomial. This polynomial is given by the
product of some irreducible factors in the resultant of~(\ref{horn})
(see Section~\ref{resultantsec}).

To prove that the function~$y(x)$ is rational it suffices to show
that~$y(x)$ depends rationally on any given variable~$x_i,$ the other
variables being fixed. Let us first consider the case when the Newton
polytope $\mathcal{N}=\mathcal{N}_{g}$ of the polynomial~$g(x)$
has zero $n$-dimensional volume (for examples of such rational
hypergeometric functions see Example~\ref{example1} in
Section~\ref{supportsec}).
Let $T\subset\R^n$ denote the minimal linear subspace whose translation
contains the polytope~$\mathcal{N}.$ Choose a basis
$u_{1},\ldots,u_{n}\in\Z^n$ of the lattice~$\Z^n$ such that
$u_{1},\ldots,u_{m}$ is a basis of the sublattice $T\cap\Z^n.$
Let us introduce new variables
$\xi_i = x^{u_i}= x_{1}^{u_{i1}} \ldots x_{n}^{u_{in}},$
$i=1,\ldots,n.$
It suffices to show that the function~$y(x(\xi))$ depends
rationally on the variables $\xi_1,\ldots,\xi_n.$

By the construction the polynomial~$g(\xi)$ is given by the product of
a monomial and another polynomial which only depends on the variables
$\xi_1,\ldots,\xi_m.$  The Newton polytope of~$g(\xi)$ has nonzero
$m$-dimensional
volume. It follows by the two-sided Abel lemma that the cone of the
support of any Laurent series $\sum_{s\in\Z^n}\varphi(s)\xi^s$
representing the function~$y(\xi)$ is contained in the linear subspace
$s_{m+1}=\ldots=s_{n}=0.$ Hence~$y(\xi)$ depends polynomially on the
variables $\xi_{m+1},\ldots,\xi_{n}.$ Let $\xi=(\xi^{'},\xi^{''}),$
where $\xi^{'}=(\xi_{1},\ldots,\xi_{m}),$ 
$\xi^{''}=(\xi_{m+1},\ldots,\xi_{n}).$
With these notations the function~$y(\xi)$ can be written in the form
$y(\xi)=\sum_{\alpha\in W} a_{\alpha}{\xi^{''}}^{\alpha} y_{\alpha}(\xi^{'}),$
where~$W$ is a finite subset of the lattice~$\Z^{n-m},$
$y_{\alpha}(\xi^{'})$ is a meromorphic function depending on the
variables $\xi_{1},\ldots,\xi_{m}$ only and $a_{\alpha}\in\C.$
We will prove that $y_{\alpha}(\xi^{'})$ is a hypergeometric function
for any $\alpha\in W.$ Since the Newton polytope of~$g(\xi)$ has nonzero
$m$-dimensional volume, this will show that it suffices to consider the
case when the Newton polytope of the polynomial defining the singular set
of a meromorphic hypergeometric function has the maximal possible dimension.

Let~$E_{i}^{\lambda_i}$ denote the operator which increases the $i$th
argument of a function depending on~$n$ variables by~$\lambda_i,$ i.e.,
$E_{i}^{\lambda_i}f(x)=f(x+\lambda_{i}e_{i}).$ For $\lambda\in\R^n$
we denote the composition of the operators
$E_{1}^{\lambda_1},\ldots,E_{n}^{\lambda_n}$ by~$E^{\lambda},$
that is $E^{\lambda}f(x) = f(x_1+\lambda_1,\ldots,x_n+\lambda_n).$ Since
the commutator $[\theta_i,x_{j}^{\lambda_j}]$ equals
$\delta_{ij}\lambda_{j}x_{j}^{\lambda_j},$ it follows that for any
polynomial~$P$ in~$n$ variables and any $\lambda\in\Z^n$
\begin{equation}
P(\theta)x^{\lambda} = x^{\lambda}(E^{\lambda}P)(\theta).
\label{change}
\end{equation}
By the definition the function~$y(x)$ is hypergeometric and hence
satisfies the Horn system~(\ref{horn}). Using the relation~(\ref{change})
and the $i$th equation of~(\ref{horn}) we compute
$$
x_{i}^2 (E_{i}^{1}P_i)(\theta)P_i(\theta)y(x)=
{(x_iP_i(\theta))}^2 y(x)= x_iP_i(\theta)Q_i(\theta)y(x)=
$$
$$
(E_{i}^{-1}Q_i)(\theta) x_iP_i(\theta)y(x)=
(E_{i}^{-1}Q_i)(\theta) Q_i(\theta)y(x).
$$
Repeating this argument~$\lambda_i$ times we arrive at the formula
\begin{equation}
x_{i}^{\lambda_i}
\left(
\prod_{j=0}^{\lambda_i -1}
(E_{i}^{j}P_i)(\theta)
\right) y(x) =
\left(
\prod_{j=0}^{\lambda_i -1}
(E_{i}^{-j}Q_i)(\theta)
\right) y(x),
\label{better}
\end{equation}
which holds for any $\lambda_i\in\N.$ For $u_{ki}\geq 0$ define
polynomials $\rho_{ki}(s)=\prod_{j=0}^{u_{ki}-1} E_{i}^{j}P_i (s)$
and $\tau_{ki}(s)=\prod_{j=0}^{u_{ki}-1} E_{i}^{-j}Q_i(s)$
(by the definition the empty product equals~$1$). For $u_{ki} < 0$
define polynomials
$\rho_{ki}(s)=\prod_{j=0}^{-u_{ki}-1} E_{i}^{-j}Q_i(s),$
$\tau_{ki}(s)=\prod_{j=0}^{-u_{ki}-1} E_{i}^{j} P_i(s).$
It follows from~(\ref{better}) that for any $k=1,\ldots,n$
\begin{equation}
x_{i}^{u_{ki}}\rho_{ki}(\theta)y(x)= \tau_{ki}(\theta)y(x), \,\,\,
i=1,\ldots,n.
\label{good}
\end{equation}
Composing the operators in the equations~(\ref{good}) in the same way
as we did before in order to obtain the formula~(\ref{better}), we
arrive at the system of equations
\begin{equation}
x^{u_k}
\left( \prod_{j=1}^{n}
\left(
\prod_{l=j+1}^n E_{l}^{u_{kl}}
\!\right) \rho_{kj}(\theta)
\!\right) y(x) =
\left( \prod_{j=1}^n
\left(
\prod_{l=1}^{j-1} E_{l}^{-u_{kl}}
\!\right) \tau_{kj}(\theta)
\!\right) y(x), \,\,\, k=1,\ldots,n.
\label{best}
\end{equation}
For instance,
$$
x_{1}^{u_{k1}}x_{2}^{u_{k2}}
(E_{2}^{u_{k2}} \rho_{k1})(\theta)\rho_{k2}(\theta) y(x) =
{\rm (by~(\ref{change}))}
=x_{1}^{u_{k1}}\rho_{k1}(\theta) x_{2}^{u_{k2}} \rho_{k2}(\theta) y(x)=
$$
$$
{\rm (by~the~2nd~equation~in~(\ref{good}))}
= x_{1}^{u_{k1}}\rho_{k1}(\theta)\tau_{k2}(\theta) y(x)=
{\rm (by~(\ref{change}))}=
$$
$$
(E_{1}^{-u_{k1}}\tau_{k2})(\theta) x_{1}^{u_{k1}}\rho_{k1}(\theta) y(x)=
{\rm (by~the~1st~equation~in~(\ref{good}))}=
$$
$$
(E_{1}^{-u_{k1}}\tau_{k2})(\theta)\tau_{k1}(\theta)y(x).
$$
Each equation in~(\ref{best}) is obtained by repeating this
argument~$n$ times.

Making in~(\ref{best}) the change of variables $\xi_i=x^{u_i}$ and using
the equality
$\theta_i= x_i\frac{\partial}{\partial x_i} =
u_{1i}\xi_1 \frac{\partial}{\partial \xi_1} + \ldots +
u_{ni}\xi_n \frac{\partial}{\partial \xi_n},$
we conclude that~$y(\xi)$ is a solution to the system of equations
\begin{equation}
\xi_i \rho^{(i)} (\theta_\xi) y(\xi) = \tau^{(i)}(\theta_\xi) y(\xi),
\,\,\, i=1,\ldots,n,
\label{myhorn}
\end{equation}
where $\theta_{\xi}=\left(\xi_1\frac{\partial}{\partial \xi_1},\ldots,
\xi_n\frac{\partial}{\partial \xi_n}  \right),$
$U$ is the matrix with the rows $u_1,\ldots,u_n$ and
$$
\rho^{(i)}(s) = \prod_{j=1}^{n}
\left( \prod_{l=j+1}^{n} E_{l}^{u_{kl}}  \right) \rho_{kj}((U^T)^{-1}s),
$$
$$
\tau^{(i)}(s) = \prod_{j=1}^{n}
\left( \prod_{l=1}^{j-1} E_{l}^{-u_{kl}}\right)  \tau_{kj}((U^T)^{-1}s).
$$
Since $y(\xi) = \sum_{\alpha\in W} a_{\alpha} {\xi^{''}}^{\alpha}
y_{\alpha}(\xi^{'}),$ it follows from the first~$m$ equations of the
system~(\ref{myhorn}) that
$$
\left( \xi_i \rho^{(i)}(\theta_\xi) - \tau^{(i)}(\theta_\xi) \right) y(\xi)=
\sum_{\alpha\in W} a_{\alpha} {\xi^{''}}^{\alpha}
\left(
\left(
\xi_i \rho^{(i)}(\theta_\xi) - \tau^{(i)}(\theta_\xi)
\right) y_{\alpha}(\xi^{'})
\right) = 0
$$
for $i=1,\ldots,m.$ Since~$y_{\alpha}(\xi^{'})$ does not depend on
$\xi_{m+1},\ldots,\xi_{n},$ it follows that for any $\alpha\in W$
\begin{equation}
\xi_i \rho^{(i)} (\theta_{\xi}^{'}) y_{\alpha}(\xi^{'}) =
\tau^{(i)}(\theta_{\xi}^{'}) y_{\alpha}(\xi^{'}), \,\,\,
i=1,\ldots,m.
\label{finalhorn}
\end{equation}
Here $\theta_{\xi}^{'}=\left(\xi_1\frac{\partial}{\partial \xi_1},\ldots,
\xi_m\frac{\partial}{\partial \xi_m},0,\ldots,0  \right).$
The system~(\ref{finalhorn}) is a Horn system in~$m$ variables. Thus the
function~$y_{\alpha}(x)$ is hypergeometric for any $\alpha\in W.$ By the
assumption the Newton polytope of the polynomial defining the singularity
of the meromorphic function~$y_{\alpha}(\xi^{'})$ has dimension~$m.$
To prove that the original function~$y(x)$ is rational it suffices to
show that~$y_{\alpha}(\xi^{'})$ depends rationally on $\xi_1,\ldots,\xi_m$
for any $\alpha\in W.$ Thus it is sufficient to prove the theorem in the
case when the Newton polytope of the polynomial which defines the singular
set of a given meromorphic hypergeometric function has the maximal
possible dimension.

Suppose now that ${\rm dim}\, \mathcal{N}=n.$ Let~$C_{v}^{\vee}$ be
the cone which is dual to~$\mathcal{N}$ at the point~$v.$ By the
remark after Theorem~B to each vertex~$v$ of the
polytope~$\mathcal{N}$ one can associate a connected component of
the amoeba complement~$^c\!\mathcal{A}_{g}.$ This component is the
image of the domain of convergence of some Laurent series~$L_{v}$
for the function $y(x)=h(x)/g(x)$ under the mapping~${\rm Log}.$
It contains some translation $w_{v}+C_{v}^{\vee}$ of the
cone~$C_{v}^{\vee}.$ By the two-sided Abel lemma the cone of the
support of the series~$L_{v}$ coincides with the cone
$-{(C_{v}^{\vee})}^{\vee}=-C_{v}.$ The family of the
cones~$\{ C_{v}^{\vee} \}_{v\in{\rm vert}(\mathcal{N})}$ coincides
with the set of all maximal cones of the dual
fan~$\Sigma_{\mathcal{N}}$ of the polytope~$\mathcal{N}.$ Since for
any polytope its dual fan is complete, it follows that the toric
variety~$\X_{\Sigma_{\mathcal{N}}}$ associated with the
fan~$\Sigma_{\mathcal{N}}$ is compact (see \S~2.4 in~\cite{Fulton}).
This variety can be covered by the affine toric
varieties~$\{ U_{C_{v}^{\vee}}\}_{v\in{\rm vert}(\mathcal{N})}.$

It is known that the monomials $\{x^{\alpha} : \alpha\in -C_{v}\}$
are holomorphic in~$U_{C_{v}^{\vee}}$ (see \S~1.3 in~\cite{Fulton}).
Since the cone of the support of the series~$L_{v}$ coincides
with~$-C_{v},$ it follows that for some~$w_{v}\in\Z^n$ the
series~$x^{w_{v}}L_{v}$ contains only those monomials which are
holomorphic in~$U_{C_{v}^{\vee}}.$
Thus~$x^{w_{v}}y(x)$ is holomorphic in~$U_{C_{v}^{\vee}}$ for all
$v\in{\rm vert}(\mathcal{N}).$ This shows that the restriction
of~$y(x)$ to any line~$x_{j}={\rm const}$ has polynomial growth
in~$\C^{*}$ and hence is rational. It is well-known that a function
which is rational in each variable depends rationally on all of the
variables. This completes the proof.~\hfill~$\square$
\end{proof}

Thanks to Theorem~\ref{merorat} we do not need to make any difference
between meromorphic and rational nonconfluent hypergeometric functions.
 From now on we formulate all the results using the term ``rational''.
\begin{remark}
\rm
Let~$f(x)$ be a rational function in~$n$ variables with singularities
along an algebraic hypersurface~$V\subset\C^n$ and let~$\mathcal{A}$ be
the image of~$V$ under the mapping~${\rm Log}.$ By Theorem~A
the connected components of the amoeba complement~$^c\!\mathcal{A}$ are
in bijective correspondence with the Laurent series expansions (with
the center at the origin) of~$f(x).$ For a multi-valued analytic
function~$F(x)$ with singularities on the same variety~$V$ this
correspondence is in general not one-to-one. It may happen that some
of the connected components of~$^c\!\mathcal{A}$ do not correspond to any
expansion of~$F(x)$ since there is no holomorphic branch of~$F(x)$ on the
pull-back of this component. It is also possible that several
connected components of~$^c\!\mathcal{A}$ correspond to a single series
expansion of~$F(x).$ (For instance, let~$x\in\C$ and consider the function
$F(x)=\sqrt{\sqrt{x+2} + \sqrt{3}}.$ There exists a holomorphic branch
of~$F(x)$ in the disk $\{ |x|<2 \}$ although~$x=1$ is a branching point.
A similar situation in the two-dimensional case is described in
Example~\ref{Cardano}.) However, with each series expansion of~$F(x)$
centered at the origin one can associate at least one connected
component of~$^c\!\mathcal{A}.$
\label{nonmero}
\end{remark}


\section{Rational solutions to the Horn system
\label{classificationsec}}

Typically a hypergeometric function is a multi-valued analytic
function with singularities along an algebraic hypersurface
(see Section~\ref{resultantsec}). In this section we give a
necessary condition for a hypergeometric series to represent a
germ of a rational function. This allows one to give an explicit
description of the class of rational solutions to~(\ref{horn})
in the case when $Q_{i}(s)=\prod_{k=0}^{p_{i}-1}(s_{i} + k/p_{i})$
for some positive integers~$p_{i},$
each linear factor of~$P_{i}(s)$ depends
on all the variables and the resultant of~(\ref{horn}) is
irreducible. We prove that any such rational hypergeometric
function is contiguous to the Bergman kernel
$K_{p_{1},\ldots,p_{n}}$ for some $p_{1},\ldots,p_{n}$
(Proposition~\ref{contigprop}).

Recall that $B_1, \ldots ,B_d$ are defined to be the duals to the maximal
elements (with respect to inclusion) of the finite family $\{-C_{I}\}$
of strongly convex polyhedral cones. Here~$C_{I}$ is the recession
cone of the convex hull of the support of the hypergeometric
series~(\ref{convergent}).
Let $X_1, \ldots ,X_N$ denote the recession cones of the connected components
of the amoeba complement~$^c\!\mathcal{A}_{R(x)}$ of the resultant
of~(\ref{horn}). These recession cones are well-defined since by
Theorem~A the connected components of the amoeba complement
are convex. The following theorem describes the structure of the
amoeba~$\mathcal{A}_{R(x)}.$

\begin{theorem}
Suppose that a nonconfluent Horn system possesses a rational solution
with the poles on the zero set of its resultant~$R(x)$. Then the fan of
this Horn system is well-defined and dual to the Newton polytope of~$R(x).$
\label{amoebathm}
\end{theorem}

\begin{proof}
Since there exists a rational solution
to~(\ref{horn}) with the poles on the zero set of its resultant~$R(x)$
it follows by Theorems~B and~\ref{rathypmin}
that the cone~$X_i$ has nonempty interior for any $i=1, \ldots ,N.$
Thus by Theorem~B the cones~$\{X_i\}_{i=1}^{N}$ can
be identified with the maximal cones of the fan which is dual to the
Newton polytope of~$R(x).$

It suffices to show that the family~$\{B_i\}_{i=1}^{d}$ consists of
the same elements as the family~$\{X_i\}_{i=1}^{N}.$
As we have already mentioned in Section~\ref{fansec} the nonconfluency
condition~(\ref{nonconfluency}) for the Horn system~(\ref{horn}) implies
that $\bigcup_{j=1}^{d} B_j = \R^n.$
Hence for any $i=1, \ldots ,N$ there exists $k_i\in\{1, \ldots ,d\}$ such that
${\rm int }(X_i\cap B_{k_i})\neq\emptyset.$
Let~$L_i$ denote a series solution to~(\ref{horn}) whose support~$S_i$
defines the cone~$B_{k_i}$ in the sense that $B_{k_i}=-C_{S_i}^{\vee}.$
Here~$C_{S_i}$ is the cone of~$S_i$ (see Section~\ref{fansec}).
Let~$\tilde{L}_i$ denote the series expansion of the rational solution
to~(\ref{horn}) such that the recession cone of the image of its domain
of convergence under the mapping ${\rm Log}$ is~$X_i.$ Since
${\rm int }(X_i\cap B_{k_i})\neq\emptyset$ it follows that the series
$L+\tilde{L}_i$ has a nonempty domain of convergence~$\Omega_i.$
By the two-sided Abel lemma the cone of the convex set
${\rm Log}\, \Omega_i$ is $X_i\cap B_{k_i}.$

Any Puiseux series solution to~(\ref{horn}) whose domain
of convergence lies entirely in the pre-image of a connected
component of the amoeba complement~$^c\!\mathcal{A}_{R(x)}$ with
respect to the mapping ${\rm Log}$ converges on the whole of this
pre-image. Using the two-sided Abel lemma we conclude that~$B_{k_i}$
cannot be a proper subset of~$X_i.$
Thus either $X_i=B_{k_i}$ or~$B_{k_i}^{\vee}$ is a proper
subset of ${(X_i \cap B_{k_i})}^{\vee}.$ The latter is impossible
due to the assumption that~$B_{k_i}^{\vee}$ is a maximal element
in the family of the cones of the supports of series solutions
to~(\ref{horn}). Hence $X_i=B_{k_i}$ for any $i=1, \ldots ,N.$
Since the cones~$\{X_i\}_{i=1}^{N}$ are the maximal cones of a complete
fan, it follows that $d=N$ and thus we can identify the families of
the cones~$\{X_i\}_{i=1}^{N}$ and~$\{B_i\}_{i=1}^{d}.$
The proof is complete.~\hfill~$\square$
\end{proof}
\vskip0.15cm

The conditions in Theorem~\ref{amoebathm} are sufficient for the fan
of a Horn system to be dual to the Newton polytope of its resultant,
but they are not necessary.
For instance, the fan of the system~(\ref{cardanohorn})
in Example~\ref{Cardano} below is dual to the Newton polytope of its
resultant though the system~(\ref{cardanohorn}) has no nonzero
rational solutions.
Yet, the remark in the very end of Section~\ref{fansec} shows that the
conclusion of Theorem~\ref{amoebathm} does not hold in arbitrary case.

\begin{corollary}
If a Horn system possesses a rational solution with the poles on the
zero set of its resultant then the number of $0$-admissible sets
associated with this system cannot be smaller than the number of the
maximal cones in its fan.
\label{merocol}
\end{corollary}

\begin{proof}
By Theorem~\ref{amoebathm} the fan of the Horn system is well-defined.
Let~$y(x)$ be a rational solution to~(\ref{horn}) with the poles on
the zero set of the resultant~$R(x)$ of~(\ref{horn}).
By Theorem~A the number of Laurent series expansions
of~$y(x)$ with the center at the origin equals the number of connected
components of the set~$^c\!\mathcal{A}_{R}.$
By Theorem~\ref{rathypmin} the amoeba of~$R(x)$ is solid and hence
by Theorem~\ref{amoebathm} there exists a one-to-one correspondence
between the connected components of~$^c\!\mathcal{A}_{R}$ and the
maximal cones of the fan of the system~(\ref{horn}).
Since any expansion of~$y(x)$ is supported in a $0$-admissible set it
follows that the number of such sets cannot be smaller than the number
of the maximal cones in the fan of the Horn system.
This completes the proof of the corollary.~\hfill~$\square$
\end{proof}
\vskip0.15cm

As we have seen in Section~\ref{resultantsec} a solution to the Horn
system~(\ref{horn}) can only be singular on the set on which the
resultant~$R(x)$ of~(\ref{horn}) vanishes. Typically~$R(x)$ is
divisible by some monomial~$x^a,$ $a\in \N^n.$ We denote the
quotient~$R(x)/x^a$ (with the maximal possible~$|a|=a_{1}+ \ldots +a_{n}$)
by~$r(x)$ and call it {\it the essential resultant} of the
system~(\ref{horn}). The reason for introducing this terminology is the
fact that a Laurent monomial has unique Laurent series development with
the center at the origin. Therefore such a monomial is an unessential
factor as long as one is concerned with the problem of computing the
number of connected components of the amoeba complement of a mapping.

The case when the polynomial~$Q_{i}(s)$ depends only on~$s_{i}$ for all
$i=1,\ldots,n$ is particularly important.
Under this assumption
it is possible to compute the dimension of the space of holomorphic
solutions to the Horn system~(\ref{horn}) explicitly and construct
a basis in this space if the parameters of the system are
sufficiently general~\cite{Sadykov}. (Theorem~9 in~\cite{Sadykov}
assumes that ${\rm deg\,} Q_{i} > {\rm deg\,}P_{i},$ $i=1,\ldots,n,$
which is not the case if the nonconfluency relation~(\ref{nonconfluency})
holds. Yet, by the lemma in~\S~1.4 of~\cite{GGR} each of the basis series
which were constructed in~\S~3 of~\cite{Sadykov} converges in some
neighbourhood of the origin if the original Horn system is nonconfluent.
The multi-valued analytic functions determined by these series give a
global basis in the space of holomorphic solutions to~(\ref{horn}).)
Recall that two Ore-Sato coefficients (and the corresponding
hypergeometric series) are called {\it contiguous} if their quotient
can be reduced to the product of a rational function and an exponential
term $\tilde{t}_{1}^{s_{1}}\ldots \tilde{t}_{n}^{s_{n}}.$
The next proposition provides an explicit description of the class
of rational solutions to such systems of hypergeometric type
under some additional assumptions on the parameters.

\begin{proposition}
Suppose that the nonconfluent Ore-Sato coefficient
$$
\psi(s) = t_{1}^{s_{1}}\ldots t_{n}^{s_{s}}
\frac
{\prod_{i=1}^{p} \Gamma(\langle A_{i},s \rangle - c_{i} )}
{\prod_{j=1}^{n}
\Gamma(p_{j}(s_{j} + 1))}
$$
defines the Horn system~(\ref{horn}) with the irreducible essential
resultant~$r(x)$ and satisfies the conditions $A_{ij}>0,$
$i=1,\ldots,p,$ $j=1,\ldots,n.$
Let $y(x) = \sum_{s\in\N^{n}} \psi(s) x^{s}$ and
let~$A$ be the matrix with the rows $A_{1},\ldots,A_{p}.$
If ${\rm rank\,} A >1$ then the series~$y(x)$ cannot define a
rational function. (We disregard exceptional values of the
parameters of~$\psi(s)$ for which~$y(x)$ reduces to a linear
combination of hypergeometric series in fewer variables.)
If ${\rm rank\,} A = 1$ and~$y(x)$ is rational then it is contiguous
to the series~(\ref{bergman}) converging to the Bergman
kernel~$K_{p_{1},\ldots,p_{n}}(x).$
\label{contigprop}
\end{proposition}

\begin{proof}
Suppose that ${\rm rank\,} A > 1$ and~$y(x)$ is a rational function.
We may without loss of generality assume that
$A_{11}A_{22}-A_{12}A_{21}\neq 0.$
For each $m=1,\ldots,p$ consider the Ore-Sato coefficient
$$
\chi_{m}(s) =
\frac
{\prod_{i=1}^{m} \Gamma(\langle A_{i},s \rangle - c_{i})}
{\prod_{j=1}^{n}
\Gamma(p_{j}(s_{j} + 1))}.
$$
Each of these coefficients defines a system of
differential equations of the Horn type (see Remark~\ref{definecoeff}).
Let $B_{m1},\ldots,B_{md_{m}}$ be the maximal elements in the family
of the cones of the admissible sets associated with the system defined
by~$\chi_{m}(s)$ (see Section~\ref{fansec}).
Arguing as in the proof of Proposition~\ref{minimalprop} we conclude
that~$d_{1}=n+1.$ Let~$\tilde{A}$ be the matrix with the rows
$A_{1},A_{2},e_{3},\ldots,e_{n},$
$\tilde{c}=(c_{1},c_{2},0,\ldots,0)\in\C^n$ and define~$\gamma$ to be
the solution to the system of linear equations $\tilde{A}s =\tilde{c}.$
The set $\{s\in\Z^n + \gamma : \tilde{A}s \geq 0\}$ satisfies the
conditions in Proposition~\ref{suppthm} if the parameters
$c_{1},\ldots,c_{p}$ are generic. This yields $d_{2}\geq n+2.$
By Remark~\ref{refinerem}
$d_{i}\leq d_{j}$ for~$i\leq j.$ Since $\chi_{p}(s)=\psi(s)$
it follows by Theorem~\ref{amoebathm} that the number of connected
components of the amoeba complement~$^c\!\mathcal{A}_{r(x)}$ at least
equals~$n+2.$ By the assumption the series~$y(x)$ represents a germ of
a rational function. Since~$r(x)$ is irreducible, the function~$y(x)$
must be singular on the whole of the hypersurface $\{r(x)=0\}.$
Thus it follows from Theorem~A that the number of
Laurent series developments (centered at the origin) of this rational
function at least equals~$n+2.$ Yet, the
condition~$A_{ij}>0$ and the conditions~(\ref{firstcond}),(\ref{secondcond})
in Proposition~\ref{suppthm} imply that the number of $0$-admissible
subsets associated with the Horn system defined by the Ore-Sato
coefficient~$\psi(s)$ cannot exceed~$n+1.$ This contradicts the
conclusion of Corollary~\ref{merocol} and shows that the function~$y(x)$
cannot be rational unless ${\rm rank\,} A = 1.$

Suppose now that ${\rm rank\,} A = 1$ and that the series~$y(x)$
converges to a rational function. Let
$\delta = {\rm GCD\,} (p_{1},\ldots,p_{n}),$
$\tilde{p}_{i} = p_{i}/\delta,$ $i=1,\ldots,n.$
It follows from the nonconfluency condition
$\sum_{i=1}^{p}A_{i}=(p_{1},\ldots,p_{n})$
and the Gauss multiplication formula for the $\Gamma$-function
that~$\psi(s)$ is contiguous to
$\tilde{\psi}(s) = \prod_{l=0}^{\delta-1}
\Gamma(\tilde{p}_{1}s_{1}+\ldots+\tilde{p}_{n}s_{n} + a_{l})/
\prod_{j=1}^{n} \Gamma(p_{j}(s_{j} + 1)).$
Here $a_{0},\ldots,a_{\delta-1}\in\C$ are some constants. Moreover the
quotient~$\psi(s)/\tilde{\psi}(s)$ is given by an exponential
term~$\tilde{t}_{1}^{s_{1}}\ldots \tilde{t}_{n}^{s_{n}}$ and hence the
series $\tilde{y}(x) = \sum_{s\in\N^{n}} \tilde{\psi}(s) x^s$
converges to a rational function. By the assumption $\tilde{p}_{i}\neq 0$
for any $i=1, \ldots ,n.$
The restriction of~$\tilde{y}(x)$ to the complex line
$x_{1}= \ldots [i] \ldots =x_{n}=0$ is a rational function
(here~$[i]$ is the sign of omission).
Let $\tilde{\psi}_{i}(s_{i})=\tilde{\psi}(0, \ldots ,s_{i}, \ldots ,0)$
($s_{i}$~in the $i$th position).
Using once again the Gauss multiplication formula we conclude that
the series
$$
\sum_{s_{i}=0}^{\infty}
\frac
{\prod_{l=0}^{\delta-1} \prod_{j=0}^{\tilde{p}_{i}-1}
\Gamma \left( s_{i} + \frac{a_{l}+j}{\tilde{p}_{i}} \right)}
{\prod_{k=0}^{p_{i}-1} \Gamma \left( s_{i} + \frac{k}{p_{i}}\right)}
x_{i}^{s_{i}}
$$
represents a rational function. A criterion for a power series in one
variable to converge to a rational function (see Theorem~4.1.1
in~\cite{Stanley}) implies that
for any $l=0, \ldots ,\delta -1,$ $j=0, \ldots ,\tilde{p}_{i}-1$ there exists
$k\in \{0, \ldots ,p_{i}-1\}$ such that
$(a_{l}+j)/\tilde{p}_{i} - k/p_{i}\in\N.$ Hence for any
$l=0, \ldots ,\delta -1$ one can find $k\in \{0, \ldots ,p_{i}-1\}$ such that
$a_{l} - k/\delta\in\Z.$
Thus~$\psi(s)$ is contiguous to the Ore-Sato coefficient
$$
\frac
{\prod_{l=0}^{\delta-1} \Gamma
\left(
\tilde{p}_{1}s_{1}+\ldots+\tilde{p}_{n}s_{n} + \frac{l}{\delta}
\right)}
{\prod_{j=1}^{n} \prod_{k=0}^{p_{j}-1}
\Gamma(s_{j} + \frac{k}{p_{j}} + 1)}.
$$
The Gauss multiplication formula shows that the latter coefficient is
contiguous to the coefficient of the series~(\ref{bergman})
which represents the Bergman kernel~$K_{p_{1},\ldots,p_{n}}.$
The proof is complete.~\hfill~$\square$
\end{proof}

\begin{remark}
\rm
There exist rational hypergeometric functions that cannot be described
in terms of the Bergman kernels of complex ellipsoidal domains. For
instance, the hypergeometric series
$$
\sum_{s\in\N_{0}^{n}}
\frac
{\Gamma(s_{1} + p(s_{2} + \ldots + s_{n} + 1))
\Gamma(s_{2} + \ldots + s_{n} + 1)}
{\Gamma(s_{1}+1)\ldots \Gamma(s_{n}+1)\Gamma(p(s_{2} + \ldots + s_{n} + 1))}
x^{s} =
\phantom{-----}
$$
$$
\phantom{---------------}
{((1-x_{1})^p - x_{2} - \ldots - x_{n})}^{-1}
$$
is not contiguous to such a kernel whenever $n\geq 3, p\geq 2.$
\label{nonbergman}
\end{remark}

Let us now consider an example. This example deals with a simplified
version of the hypergeometric series which expresses a solution~$y(x)$
to the cubic equation $y^3 + x_{1}y^2  + x_{2}y -1 = 0$ in terms of
the coefficients~$x_{1},x_{2}$
(see~\cite{Mellin},\cite{Semusheva},\cite{Sturmfels} and
Corollary~\ref{discriminant}).
\begin{example}
\rm
Consider the hypergeometric series
\begin{equation}
y(x_{1},x_{2}) = \sum_{s_{1},s_{2}\geq 0}
\frac
{\Gamma (2s_{1}+s_{2}+\alpha)\Gamma(s_{1}+2s_{2}+\beta)}
{\Gamma(3s_{1}+3) \Gamma(3s_{2}+3)} x_{1}^{s_{1}} x_{2}^{s_{2}},
\label{cardanoseries}
\end{equation}
where~$\alpha,\beta$ are arbitrary parameters such that the coefficient
of the series is well-defined and different from zero on~$\N_{0}^{2}.$
By the lemma in~\S~1.4 of~\cite{GGR} the series~(\ref{cardanoseries})
converges in some neighborhood of the origin. This series satisfies
the system of equations of hypergeometric type
\begin{equation}
\left\{
\begin{array}{clcr}
x_{1}(2\theta_{1}+\theta_{2}+\alpha)(2\theta_{1}+\theta_{2}+\alpha +1)
(\theta_{1}+2\theta_{2}+\beta) y(x) = \phantom{---} \\
\phantom{---------------} 3\theta_{1}(3\theta_{1}+1)(3\theta_{1}+2)y(x), \\
x_{2}(2\theta_{1}+\theta_{2}+\alpha)(\theta_{1}+2\theta_{2}+\beta)
(\theta_{1}+2\theta_{2}+\beta+1) y(x) = \phantom{---} \\
\phantom{---------------} 3\theta_{2}(3\theta_{2}+1)(3\theta_{2}+2)y(x).
\end{array}
\right.
\label{cardanohorn}
\end{equation}
The principal symbols of the operators in~(\ref{cardanohorn}) are
$$
\begin{array}{clcr}
H_{1}(x,z) =& x_{1}(2x_{1}z_{1}+x_{2}z_{2})^2
(x_{1}z_{1}+2x_{2}z_{2}) - 27(x_{1}z_{1})^3, \\
H_{2}(x,z) =& x_{2}(2x_{1}z_{1}+x_{2}z_{2})
(x_{1}z_{1}+2x_{2}z_{2})^2 - 27(x_{2}z_{2})^3.
\end{array}
$$
The singular locus of a solution to~(\ref{cardanohorn}) is contained
in the set on which the polynomials $H_{1}(x,z),H_{2}(x,z)$ (considered
as polynomials in~$z_{1},z_{2}$ whose coefficients depend on the
parameters~$x_{1},x_{2}$) do not form a regular sequence (see
the remark after the proof of Proposition~\ref{algebraicsing}). This
happens if and only if the resultant of $H_{1}(x,z),H_{2}(x,z)$ with
respect to~$z_{1},z_{2}$ is equal to zero. This resultant is given by
$$
R(x_{1},x_{2}) = x_{1}^{9}x_{2}^{9}
(x_{1}^{2}x_{2}^{2} + 64x_{1}^{3} - 24x_{1}^{2}x_{2} - 24x_{1}x_{2}^{2}
+ 64x_{2}^{3}
$$
\begin{equation}
- 1296x_{1}^{2} + 4698x_{1}x_{2} - 1296x_{2}^{2}
+ 8748x_{1} + 8748x_{2} - 19683).
\label{cardanoresultant}
\end{equation}
The essential resultant $r(x_{1},x_{2})=R(x_{1},x_{2})/(x_{1}x_{2})^9$
of the system~(\ref{cardanohorn}) is an irreducible polynomial.
The vectors $(2,1),(1,2)$ of the coefficients of the linear factors in
the arguments of the $\Gamma$-functions in the numerator of the
coefficient of~(\ref{cardanoseries}) are linearly independent. By
Proposition~\ref{contigprop} the series~(\ref{cardanoseries}) cannot
converge to a rational function.

\begin{minipage}{4cm}
\vskip0.2cm
\begin{picture}(80,80)
   \put(40,0){\vector(0,1){80}}
   \put(0,40){\vector(1,0){80}}
   \put(80,20){\line(-2,1){80}}
   \put(20,80){\line(1,-2){40}}

\setshadegrid span <1pt>

\hshade   40 40  70 70  40 40 /
\hshade   20 40  80 40  40 40 /
\setshadegrid span <3pt>
\hshade  40  0  40 60   0  0  /
\hshade   0  0  40 40   0 40  /
\hshade   0 40  60 40  40 40  /

\end{picture}

\noindent
{\scriptsize
{\bf Fig.~4} The maximal cones of the irreducible supports of solutions
to~(\ref{cardanohorn})
}
\end{minipage}
\hskip0.75cm
\begin{minipage}{4cm}
\hskip0.1cm
\scalebox{0.75}[0.75]

{\includegraphics*{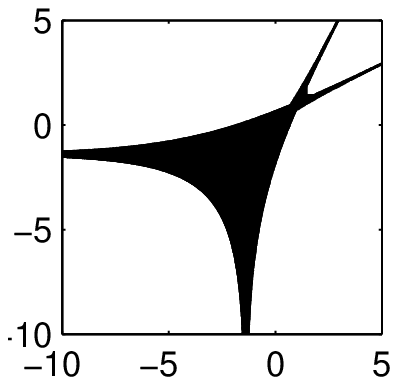}}

\hskip-3.71cm \vskip-0.72cm - \vskip0.25cm

\noindent
{\scriptsize
{\bf Fig.~5} The amoeba of the essential resultant of the Horn
system~(\ref{cardanohorn})
}
\end{minipage}
\hskip0.75cm
\begin{minipage}{4cm}
\begin{picture}(80,80)
   \put(0,0){\line(0,1){80}}
   \put(0,0){\line(1,0){80}}
   \put(80,0){\line(0,1){80}}
   \put(0,80){\line(1,0){80}}
   \put(60,0){\line(-1,2){20}}
   \put(40,40){\line(-2,1){40}}
   \put(0,0){\circle*{3}}
   \put(20,0){\circle*{3}}
   \put(40,0){\circle*{3}}
   \put(60,0){\circle*{3}}
   \put(0,20){\circle*{3}}
   \put(20,20){\circle*{3}}
   \put(40,20){\circle*{3}}
   \put(0,40){\circle*{3}}
   \put(20,40){\circle*{3}}
   \put(40,40){\circle*{3}}
   \put(0,60){\circle*{3}}
\end{picture}

\noindent
{\scriptsize
{\bf Fig.~6} The Newton polytope of the essential resultant of the
system~(\ref{cardanohorn})
}
\end{minipage}

\vskip0.2cm
The fact that the series~(\ref{cardanoseries}) cannot define a germ
of a rational function can be seen without appealing to
Proposition~\ref{contigprop} since we have the explicit
expression~(\ref{cardanoresultant}) for the resultant of the
principal symbols of the differential operators in~(\ref{cardanohorn}).
If the sum~$y(x)$ of the series~(\ref{cardanoseries})
was rational then by Theorem~\ref{rathypmin} the number of expansions
of~$y(x)$ into a Laurent series with the center at the origin would be
equal to~$4$ since the Newton polytope of the essential resultant
of~(\ref{cardanohorn}) has~$4$ vertices (see Figure~6).
However, Proposition~\ref{suppthm} shows that for any choice of the
parameters~$\alpha,\beta$ at most~$3$ of the admissible subsets can
belong to~$\Z^2$ (see Figure~4). Thus the sum of the
series~(\ref{cardanoseries}) is not a rational function.

To determine the resultant of a general Horn system is a problem of
great computational complexity. Theorem~\ref{amoebathm} and the
corollary to it allow one to describe the amoeba of the resultant of
a Horn system and draw consequences on its solvability in the class
of rational functions without performing this computation.
\label{Cardano}
\end{example}

Finally we give an example which illustrates how Theorem~\ref{rathypmin}
(or its Corollaries~\ref{adiscmin} and~\ref{discriminant}) can be
applied to the problem of constructing polyhedral decompositions of
the Newton polytopes of discriminants.

In~\cite{PR} a natural polyhedral decomposition
of the Newton polytope of a Laurent polynomial~$f$ is given.
This decomposition is determined by the piecewise linear convex
function constructed from the so-called Ronkin function~$N_{f}(t)$
which is a convex function in~$t\in\R^n.$ The function~$N_{f}$
is affine-linear on each connected component of~$^c\!\mathcal{A}_{f}.$
If such a component~$M$ corresponds to a vertex $\nu=\nu(M)$
(see Theorem~B) of the Newton polytope of~$f,$ then the Ronkin
function~$N_{f}$ is given, for $t\in M$, by
$N_{f}(t) = \log |c_\nu| + \langle t, \nu \rangle$, where $c_\nu$
denotes the coefficient of $x^\nu$ in $f$. (See Theorem~2 in \cite{PR}
for an explanation of this.)

\begin{example}
\rm
Consider the quartic equation
\begin{equation}
y^{4} + x_1 y^{3} + x_2 y^2 + x_3 y - 1 = 0.
\label{quartic}
\end{equation}
The discriminant of~(\ref{quartic}) is given by the polynomial
$$
x_{1}^{2} x_{2}^{2} x_{3}^{2} - 4 x_{1}^{3} x_{3}^{3}
+ 4 x_{1}^{2} x_{2}^{3} - 4 x_{2}^{3} x_{3}^{2}
- 18 x_{1}^{3} x_{2} x_{3} + 18 x_{1} x_{2} x_{3}^{3}
- 27 x_{1}^{4} - 16 x_{2}^{4} - 27 x_{3}^{4} +
$$
\begin{equation}
80 x_{1} x_{2}^{2} x_{3} + 6 x_{1}^{2} x_{3}^{2}
+ 144 x_{1}^{2} x_{2} - 144 x_{2} x_{3}^{2}
- 192 x_{1} x_{3} - 128 x_{2}^{2} - 256.
\label{quarticdiscriminant}
\end{equation}
By Corollary~\ref{discriminant} the zero locus of the
polynomial~(\ref{quarticdiscriminant}) has a solid amoeba.
The Newton polytope of~(\ref{quarticdiscriminant})
is displayed in Figure~7.
\vskip0.4cm

\begin{center}
\begin{minipage}{4.5cm}
\vskip0.4cm
\begin{picture}(120,120)
   \put(60,60){\vector(-1,-1){50}}
     \put(17,10){$\scriptscriptstyle x_1$}
   \put(60,60){\line(1,0){8}} \put(72,60){\vector(1,0){48}}
     \put(110,52){$\scriptscriptstyle x_2$}
   \put(60,60){\line(0,1){4}} \put(60,67){\vector(0,1){53}}
     \put(64,115){$\scriptscriptstyle x_3$}
   \put(30,30){\line(4,1){52}}
   \put(82,43){\line(1,1){18}}
   \put(30,30){\line(1,4){10}}
   \put(100,60){\line(-1,2){10}}
   \put(90,80){\line(-3,2){30}}
   \put(60,100){\line(-2,-3){20}}
   \put(40,70){\line(4,-1){28}}
   \put(68,63){\line(2,-3){13}}
   \put(68,63){\line(5,4){22}}
   \put(100,60){\circle*{3}} \put(101,66){$\scriptscriptstyle (0,4,0)$}
   \put(82,43){\circle*{3}}  \put(85,37){$\scriptscriptstyle (2,3,0)$}
   \put(30,30){\circle*{3}}  \put(33,24){$\scriptscriptstyle (4,0,0)$}
   \put(60,100){\circle*{3}} \put(63,102){$\scriptscriptstyle (0,0,4)$}
   \put(68,63){\circle*{3}}  \put(73,62){$\scriptscriptstyle (2,2,2)$}
   \put(90,80){\circle*{3}}  \put(93,83){$\scriptscriptstyle (0,3,2)$}
   \put(40,70){\circle*{3}}  \put(14,73){$\scriptscriptstyle (3,0,3)$}
   \put(60,60){\circle*{3}}
\end{picture}

\noindent
{\scriptsize
{\bf Fig.~7} The Newton polytope of the discriminant of the
equation~(\ref{quartic})
}
\end{minipage}
\end{center}

\vskip0.3cm From the solidness of the amoeba of the
discriminant~(\ref{quarticdiscriminant}) we conclude that any affine
linear part of the function~$N_{f}$ corresponds to one of the eight
vertices of the Newton polytope of~(\ref{quarticdiscriminant}).
Taking the maximum of these eight affine linear functions, we obtain
the piecewise linear convex function
$$
\max\,\bigl(\,8\log 2\,,\ 3\log 3 + 4t_1\,,\ 4\log 2 + 4t_2\,,\ 3\log3 +
4t_3\,,\ 2\log 2 + 2t_1 + 3t_2\,, \phantom{-----}
$$
\begin{equation}\label{ronkin}
\phantom{-----------} 2\log 2 + 3t_1 + 3t_3\,,\ 2\log 2 + 3t_2 + 2t_3\,,\
2t_1+2t_2+2t_3\,\bigr)\,.
\end{equation}
The set of all points $t$ at which the convex convex function
(\ref{ronkin})
is not smooth
is a two-dimensional polyhedral complex called the {\it spine} of the
amoeba, and the Legendre transform of (\ref{ronkin}) similarly gives
rise to a dual polyhedral subdivision of the polytope in Figure~7.
It deserves to be mentioned that in this example the polyhedral
decomposition of the polytope is not simplicial, for it contains
a polytope with~5 vertices, namely the convex hull of the points
$(0,4,0), (2,3,0), (3,0,3), (0,3,2), (2,2,2).$ This is because there
is a point, $t=(3\log2,4\log2,3\log2)$, at which the maximum in
(\ref{ronkin}) is attained simultaneously by the five functions
$4\log 2 + 4t_2$, $2\log 2 + 2t_1 + 3t_2$, $2\log 2 + 3t_1 + 3t_3$,
$2\log 2 + 3t_2 + 2t_3$, and $2t_1+2t_2+2t_3$.
\label{quarticex}
\end{example}


\vskip1cm

\noindent
\begin{minipage}{8cm}
{\scriptsize
DEPARTMENT OF MATHEMATICS\\
UNIVERSITY OF STOCKHOLM\\
S-10691 STOCKHOLM\\
SWEDEN
}

\vskip0.3cm

\noindent {\small {\it E-mail:} \, passare@math.su.se

}

\vskip4.7cm

\phantom{-}

\end{minipage}
\begin{minipage}{7cm}
{\scriptsize
DEPARTMENT OF MATHEMATICS\\
KRASNOYARSK STATE UNIVERSITY\\
660041 KRASNOYARSK\\
RUSSIA
}

\vskip0.3cm

\noindent {\small {\it E-mail:} \, tsikh@lan.krasu.ru

\hskip1.67cm sadykov@lan.krasu.ru
}

\vskip4.2cm

\phantom{-}

\end{minipage}

\end{document}